\documentclass{article}
\pdfoutput=1

\usepackage[narrow]{arxiv}

\usepackage[utf8]{inputenc} % allow utf-8 input
\usepackage[T1]{fontenc}    % use 8-bit T1 fonts
\usepackage[hypertexnames=false]{hyperref}       % hyperlinks
\usepackage{url}            % simple URL typesetting
\usepackage{booktabs}       % professional-quality tables
\usepackage{amsfonts}       % blackboard math symbols
\usepackage{nicefrac}       % compact symbols for 1/2, etc.
\usepackage{microtype}      % microtypography

\usepackage{lipsum}         % Can be removed after putting your text content
\usepackage{graphicx}
\usepackage{natbib}
\usepackage{doi}

\title{Nesterov Finds GRAAL: Optimal and Adaptive Gradient Method for Convex Optimization}
\renewcommand{\headeright}{}

\newif\ifuniqueAffiliation
\uniqueAffiliationtrue

\ifuniqueAffiliation % Standard variant of author block
\author{
  Ekaterina Borodich \\
  Yandex Research \\
  \texttt{borodich.ed@phystech.edu}
  \And
  Dmitry Kovalev \\
  Yandex Research \\
  \texttt{dakovalev1@gmail.com}
}
\else
\usepackage{authblk}

\author[1]{%
  Dmitry Kovalev\thanks{\texttt{dakovalev1@gmail.com}}%
}
\affil[1]{Yandex Research}
\fi
\usepackage{algorithm,algpseudocodex}
\usepackage{nicematrix}

\usepackage[textsize=scriptsize,color=orange!20,bordercolor=orange]{todonotes}
\usepackage{xcolor}
\hypersetup{
  colorlinks=true,
  linkcolor=red!75!black,
  citecolor=blue!50!black
}

\usepackage{mathtools,amsthm,amssymb}
\allowdisplaybreaks[4]

\usepackage{enumitem}
\setlist{topsep=0pt,partopsep=0pt}

\usepackage{cleveref}
\usepackage{crossreftools}
\crefname{problem}{problem}{problems}
\crefname{condition}{condition}{conditions}
\crefname{property}{property}{properties}
\creflabelformat{problem}{(#2#1#3)}

\usepackage{pifont}
\newcommand{\cmark}{{\color{green!70!black}\ding{52}}}
\newcommand{\xmark}{{\color{red}\ding{56}}}

\usepackage{annotate}
\usepackage{mythm}
\usepackage{mymath}

\newcommand{\sX}{\R^d}
\newcommand{\bg}{\operatorname{D}}
\newcommand{\ox}{\overline{x}}

\newcommand{\tx}{\tilde{x}}
\newcommand{\hx}{\hat{x}}

\begin{document}
\maketitle

\begin{abstract}
  In this paper, we focus on the problem of minimizing a continuously differentiable convex objective function, $\min_x f(x)$. Recently, \citet{malitsky2020golden,alacaoglu2023beyond} developed an adaptive first-order method, GRAAL. This algorithm computes stepsizes by estimating the local curvature of the objective function without any line search procedures or hyperparameter tuning, and attains the standard iteration complexity $\cO\brr{{L\sqn{x_0-x^*}/\epsilon}}$ of fixed-stepsize gradient descent for $L$-smooth functions. However, a natural question arises: is it possible to accelerate the convergence of GRAAL to match the optimal complexity $\cO\brr{\sqrt{L\sqn{x_0-x^*}/\epsilon}}$ of the accelerated gradient descent of \citet{nesterov1983method}? Although some attempts have been made by \citet{li2025simple,suh2025adaptive}, the ability of existing accelerated algorithms to adapt to the local curvature of the objective function is highly limited. We resolve this issue and develop GRAAL with Nesterov acceleration, which can adapt its stepsize to the local curvature at a geometric, or linear, rate just like non-accelerated GRAAL. We demonstrate the adaptive capabilities of our algorithm by proving that it achieves near-optimal iteration complexities for $L$-smooth functions, as well as under a more general $(L_0,L_1)$-smoothness assumption \citep{zhang2019gradient}.

  % We prove that our algorithm achieves the desired optimal convergence rate for Lipschitz smooth functions. Moreover, in contrast to existing methods, it does so with an arbitrary, even excessively small, initial stepsize at the cost of a logarithmic additive term in the iteration complexity.

  % In this paper, we focus on the problem of minimizing a continuously differentiable convex objective function $\min_x f(x)$. Recently, several adaptive gradient methods, including GRAAL \citep{malitsky2020golden}, have been developed. These methods estimate the local curvature of the objective function to compute stepsizes, attain the standard convergence rate $\cO(1/k)$ of fixed-stepsize gradient descent for Lipschitz-smooth functions, and do not require any line search procedures or hyperparameter tuning. However, a natural question arises: is it possible to accelerate the convergence of these algorithms to match the optimal rate $\cO(1/k^2)$ of the accelerated gradient descent of \citet{nesterov1983method}? Although some attempts have been made \citep{li2025simple}, the capabilities of the existing accelerated algorithms to adapt to the curvature of the objective function are highly limited. Consequently, we provide a positive answer to this question and develop GRAAL with Nesterov acceleration. We prove that our algorithm achieves the desired optimal convergence rate for Lipschitz smooth functions. Moreover, in contrast to existing methods, it does so with an arbitrary, even excessively small, initial stepsize at the cost of a logarithmic additive term in the iteration complexity.
\end{abstract}

\section{Introduction}\label{sec:intro}

First-order, or gradient, optimization methods are highly popular in many practical applications due to their simplicity and scalability. However, the key limitation of these methods is that they require the choice of the stepsize parameter. In this paper, we focus on developing efficient gradient methods that can adjust the stepsize at each iteration in an adaptive manner. Formally speaking, we consider the following optimization problem:
\begin{equation}\label[problem]{eq:main}
  \min_{x \in \sX} f(x),
\end{equation}
where $\sX$ is a finite-dimensional Euclidean space, and $f(x)\colon \sX \to \R$ is a convex continuously differentiable objective function. We assume that \cref{eq:main} has a solution $x^* \in \sX$.

\subsection{Gradient Methods}

The simplest and most fundamental example of first-order methods is gradient descent (GD). This algorithm performs iterations to find an approximate solution to \cref{eq:main} according to the following update rule:
\begin{equation}\label{eq:GD}
  x_{k+1} = x_k - \eta \nabla f(x_k),
\end{equation}
where $\eta > 0$ is the stepsize. Despite its simplicity, GD and its variants are widely used in practice, especially for solving large-scale problems that often appear, for instance, in machine learning.
It is well-known \citep{polyak1963gradient,nesterov2018lectures,drori2014performance} that in the case where the objective function $f(x)$ is $L$-smooth, i.e., the gradient $\nabla f(x)$ is $L$-Lipschitz, GD with the stepsize $\eta = 1/L$ achieves the following iteration complexity:
\begin{equation}\label{eq:rate_GD}
  K = \cO\brr*{{L\sqn{x_0 - x^*}/\epsilon}}
  \quad\Rightarrow\quad
  f(x_K) - f(x^*) \leq \epsilon.
\end{equation}
In addition, in his seminal work, \citet{nesterov1983method} proposed a modification of GD that implements acceleration via momentum. This accelerated gradient descent (AGD) achieves substantially improved iteration complexity:
\begin{equation}\label{eq:rate_AGD}
  K = \cO\big(\sqrt{L\sqn{x_0 - x^*}/\epsilon}\big)
  \quad\Rightarrow\quad
  f(x_K) - f(x^*) \leq \epsilon.
\end{equation}
It was shown that AGD is an optimal algorithm. That is, the complexity in \cref{eq:rate_AGD} cannot be improved by any first-order optimization method due to the lower complexity bounds of \citet{nesterov2018lectures}.

\subsection{Adaptive Methods}\label{sec:intro:adaptive}

One of the main issues with the standard GD and AGD is that they require tuning the stepsize $\eta$. In particular, they require knowledge of the gradient Lipschitz constant $L$ to achieve the iteration complexities in \cref{eq:rate_GD,eq:rate_AGD}. A typical approach to addressing this issue is to use a time-varying stepsize $\eta_k$ in \cref{eq:GD}, which is computed at each iteration according to a certain adaptive rule.

\textbf{Line search.}
The simplest way to compute the stepsize $\eta_k$ is to use line search (or backtracking), an iterative procedure that finds $\eta_k$ satisfying a certain objective function value decrement condition. It was originally proposed by \citet{goldstein1962cauchy,armijo1966minimization}, and its modern variant for GD and AGD was analyzed by \citet{nesterov2013gradient}. Unfortunately, line search makes the iterations of gradient methods more expensive as it requires the computation of the gradient $\nabla f(x)$ and/or function value $f(x)$ multiple times without making any ``progress''. Hence, it is rarely used in practice.

\textbf{AdaGrad-type methods.} An alternative approach is to use the following stepsizes $\eta_k$ in \cref{eq:GD}:
\begin{equation}\label{eq:AdaGrad}
  \eta_k = \eta\cdot\big(\tsum_{i=0}^{k}\sqn{\nabla f(x_i)}\big)^{-1/2},
\end{equation}
where $\eta > 0$ is a positive parameter.
The resulting algorithm is called AdaGrad and was originally developed by \citet{duchi2011adaptive,mcmahan2010adaptive}. It is well known that AdaGrad has the complexity of GD in \cref{eq:rate_GD} for $L$-smooth functions \citep{levy2018online}. Moreover, this algorithm is universal: it can also achieve the corresponding complexities of GD for non-smooth but Lipschitz functions, or in the case where only stochastic estimates of the gradients are available, all with a single choice of the parameter $\eta \propto \norm{x_0 - x^*}$ \citep{levy2018online,li2019convergence,orabona2023normalized}.
In addition, accelerated variants of AdaGrad are available \citep{levy2018online,cutkosky2019anytime,kavis2019unixgrad,rodomanov2024universality,kreisler2024accelerated,kovalev2025sgd}, as well as parameter-free variants \citep{cutkosky2018black,orabona2021parameter,defazio2023learning,mishchenko2023prodigy,ivgi2023dog,khaled2023dowg,kreisler2024accelerated}, which do not require tuning the parameter $\eta$. Unfortunately, AdaGrad-type methods have a significant drawback: the stepsize $\eta_k$ in \cref{eq:AdaGrad} is non-increasing. Hence, it cannot truly adapt to the local curvature of the objective function, which may limit its performance in many applications \citep{defazio2022grad}.

\textbf{Local curvature estimation.}
In this paper, we focus on a different approach to computing the stepsize $\eta_k$ by estimating the local curvature, i.e., the local gradient Lipschitz constant. To the best of our knowledge, the first such algorithm that has strong convergence guarantees, GRAAL, was proposed by \citet{malitsky2020golden}. It is a modification of GD, which uses the following stepsize rule at each iteration:
\begin{equation}\label{eq:GRAAL}
  \eta_{k+1} = \min\brf*{(1+\gamma)\eta_k, \tfrac{\nu \lambda_{k+1}^2}{\eta_{k-1}}},
\end{equation}
where $\lambda_{k+1} > 0$ is a certain finite-difference estimate of the local inverse gradient Lipschitz constant at the current iteration, and $\gamma,\nu > 0$ are positive constants. \cite{alacaoglu2023beyond} showed that GRAAL can achieve the iteration complexity in \cref{eq:rate_GD} for $L$-smooth functions. Moreover, \citet{malitsky2020adaptive} established the same result for AdGD, the vanilla GD with a stepsize rule similar to \cref{eq:GRAAL}:
\begin{equation}\label{eq:AdGD}
  \eta_{k+1} = \min\brf*{\eta_k\sqrt{1 + \tfrac{\gamma\eta_k}{\eta_{k-1}}}, \nu\lambda_{k+1}}.
\end{equation}
Overall, GRAAL and AdGD demonstrate attractive results, both theoretically and experimentally, on a range of practical optimization problems \citep{alacaoglu2023beyond,malitsky2020adaptive}.

\subsection{Main Contribution: GRAAL with Nesterov Acceleration}

Motivated by the attractive theoretical and practical results for GRAAL and AdGD, we pose the following natural research question:
\begin{question}<q>
  Is it possible to develop an algorithm that incorporates Nesterov acceleration and can truly adapt to the local curvature of the objective function, as GRAAL and AdGD do?
  % Is it possible to accelerate the convergence of such methods in the same way that \\\citet{nesterov1983method} accelerated the convergence of fixed-stepsize GD?
\end{question}
Unfortunately, to the best of our knowledge, there is no positive and comprehensive answer to this question. \citet{malitsky2020adaptive} proposed an accelerated version of AdGD, which showed strong experimental results. However, it is only a heuristic and does not have any theoretical convergence guarantees whatsoever. In addition, \citet{li2025simple} developed AC-FGM and \citet{suh2025adaptive} developed AdaNAG, which can be seen as attempts to incorporate Nesterov acceleration into GRAAL/AdGD with theoretical guarantees. However, the abilities of AC-FGM and AdaNAG to adapt to the local curvature of the objective function are highly limited, as we will discuss later.

In this paper, we provide a positive answer to \Cref{q} and make the following contributions:
\begin{enumerate}[label=\bf(\roman*)]
  \item In \Cref{sec:alg}, we develop Accelerated GRAAL (\Cref{alg}) for solving \cref{eq:main}, which incorporates Nesterov acceleration and utilizes a generalized version of the stepsize update rules in \cref{eq:GRAAL,eq:AdGD}. We also provide a theoretical convergence analysis of this algorithm.

  \item In \Cref{sec:L}, we show that \Cref{alg} achieves the optimal iteration complexity in \cref{eq:rate_AGD} for $L$-smooth functions up to additive logarithmic factors, without the requirement of hyperparameter tuning or any additional line search procedures.

  \item In \Cref{sec:L01}, we demonstrate the adaptive capabilities of \Cref{alg} by showing that it achieves the iteration complexity in \cref{eq:rate_AGD} under the more general $(L_0,L_1)$-smoothness of the objective function, up to constant additive factors that do not depend on the precision $\epsilon$.
\end{enumerate}

The important feature of \Cref{alg} is that it can adapt the stepsize $\eta_k$ to the local curvature at a geometric, or linear, rate, just like the standard non-accelerated GRAAL. In contrast, AC-FGM and AdaNAG allow only sublinear growth of the stepsize, so their adaptive abilities are insufficient. In particular, \Cref{alg} is the first adaptive algorithm that can achieve near-optimal iteration complexity for $(L_0,L_1)$-smooth functions, while there are no such results for AC-FGM or AdaNAG, to the best of our knowledge. More details are available in \Cref{sec:L:comparison,sec:L01:comparison}.

\subsection{Related Work}\label{sec:related}

\textbf{More adaptive stepsizes: Barzilai-Borwein and Polyak.} The idea of computing stepsizes using estimates of the local gradient Lipschitz constant was previously used by \citet{barzilai1988two}, who proposed the following stepsize rule:
\begin{equation}\label{eq:BB}
  \eta_{k+1} = \tfrac{\<x_{k+1} - x_{k},\nabla f(x_{k+1}) - \nabla f(x_k)>}{\sqn{\nabla f(x_{k+1}) - \nabla f(x_k)}}.
\end{equation}
Unfortunately, GD with this rule provably works only in the case where the objective function $f(x)$ is quadratic \citep{raydan1993barzilai,dai2002r}, and may not work otherwise \citep{burdakov2019stabilized}. \citet{polyak1969minimization} suggested using GD with the following stepsize rule:
\begin{equation}
  \eta_{k+1} = \tfrac{f(x_{k+1}) - f(x^*)}{\sqn{\nabla f(x_{k+1})}}.
\end{equation}
Similar to AdaGrad, GD with this rule was shown to be universal \citep{hazan2019revisiting}. However, it requires a tight estimate of the optimal objective function value $f(x^*)$, which is rarely available in practice.

\textbf{Optimization for $(L_0,L_1)$-smooth functions.} The $(L_0,L_1)$-smoothness assumption was proposed by \citet{zhang2019gradient} as a generalization and a more realistic alternative to the standard $L$-smoothness. The convergence of gradient methods under this assumption has been extensively studied in the literature \citep{zhang2020improved,chen2023generalized}. \citet{gorbunov2024methods} showed that AdGD can achieve the iteration complexity of non-accelerated GD in \cref{eq:rate_GD} up to additive constant factors without the requirement of hyperparameter tuning or line search. Additionally, several accelerated algorithms with theoretical guarantees are available \citep{li2023convex,gorbunov2024methods,vankov2024optimizing}. However, all these algorithms are non-adaptive, and only \citet{vankov2024optimizing} managed to achieve the optimal iteration complexity in \cref{eq:rate_AGD} up to additive constant factors, with the requirement of a substantially more complex small-dimensional relaxation oracle \citep{nesterov2021primal}. It is also worth mentioning the concurrent work of \citet{tyurin2025near}, who managed to achieve the complexity in \cref{eq:rate_AGD} up to additive constant factors. However, their algorithm requires the tuning of several parameters and is, therefore, also non-adaptive. Additionally, the initial version of our paper appeared online prior to the work of \citet{tyurin2025near}, with only the results from \Cref{sec:L01} missing, which we were finalizing at that time.

\section{Adaptive Gradient Method with Nesterov Acceleration}\label{sec:alg}

\subsection{Algorithm Development}\label{sec:alg:dev}

In this section, we develop Accelerated GRAAL for solving \cref{eq:main}. Below, we briefly describe the key ideas used in the development of the algorithm. After assembling these ideas, we obtain the resulting \Cref{alg}.

\textbf{Local curvature estimator.} As discussed in \Cref{sec:intro:adaptive}, the stepsize rule in \cref{eq:GRAAL} used in GRAAL requires an estimate of the inverse local gradient Lipschitz constant $\lambda$. When the objective function is convex, given two points $x,z \in \sX$, we can consider the following two options for computing $\lambda$:
\begin{equation}\label{eq:options}
  \text{Option I:}\quad\lambda = \tfrac{\norm{x-z}}{\norm{\nabla f(x) - \nabla f(z)}},\qquad
  \text{Option II:}\quad\lambda = \tfrac{2\bg_f(x,z)}{\sqn{\nabla f(x) - \nabla f(z)}}.
\end{equation}
GRAAL was originally developed for solving monotone variational inequalities (VI). Hence, it uses Option~I, which works for this more general problem class, with the gradient $\nabla f(x)$ replaced by the monotone operator. However, it turns out that Option~II can better exploit the properties of the objective function $f(x)$. Hence, we use Option~II and, for convenience, define $\Lambda(x;z)$ as follows:
\begin{equation}\label{eq:Lambda}
  \Lambda(x;z) =
  \begin{cases}
    \frac{2\bg_f(x;z)}{\sqn{\nabla f(x) - \nabla f(z)}} & \nabla f(x) \neq \nabla f(z)\\
    +\infty & \nabla f(x) = \nabla f(z)
  \end{cases}.
\end{equation}
It is worth noting that \citet{li2025simple,suh2025adaptive} also used Option~II in AC-FGM and AdaNAG, respectively.

\textbf{Nesterov acceleration.} To incorporate Nesterov acceleration in \Cref{alg}, we use its recent interpretation by \citet{kovalev2024linear}. The idea is that at the iteration $k$ of a gradient method, we replace the objective function $f(x)$ with the function $f_k(x)\colon \sX \to \R$, defined as follows:
\begin{equation}\label{eq:f_k}
  f_k(x) = \alpha_k^{-1} \cdot f(\alpha_k x + (1-\alpha_k)\ox_k),
  \quad\text{where}\;\;
  \alpha_k \in (0,1],\;\;
  \ox_k \in \sX.
\end{equation}
In the case where GD is used, that is, $x_{k+1} = x_k - \eta_k \nabla f_k(x_k)$, we can choose $\alpha_k = 2/(k+2)$ and $\ox_{k+1} = \alpha_k x_{k+1}+ (1-\alpha_k)\ox_k$, which gives the STM algorithm \citep{gasnikov2016universal}, a variant of AGD. However, we use the definition in \cref{eq:f_k} as it substantially simplifies the development of \Cref{alg}.

\textbf{GRAAL extrapolation.} We use the extrapolation step of GRAAL in combination with the interpretation of Nesterov acceleration above. It can be summarized as follows:
\begin{equation}\label{eq:extrapolation}
  x_{k+1} = x_k - \eta_k \nabla f_k(\hx_k),\qquad
  \hx_{k+1} = x_{k+1} + \theta(x_{k+1} - x_k),
\end{equation}
where $\theta > 0$ is the extrapolation parameter. GRAAL uses extrapolation for two reasons. First, the vanilla gradient method does not work for VI as it diverges even on simple bilinear min-max problems \citep{daskalakis2017training}. Hence, the extragradient method \citep{korpelevich1976extragradient,mishchenko2020revisiting} or methods with extrapolation \citep{daskalakis2017training,malitsky2020forward,kovalev2022accelerated} are typically used. Second, and more importantly, to the best of our knowledge, the particular type of extrapolation used by GRAAL plays a key role in its adaptive capabilities. In particular, it is an open question whether our results can be obtained with a different baseline algorithm.

\textbf{Problem: choosing $\alpha_k$.} Although it may seem that the tools described above are already enough to obtain \Cref{alg}, one issue remains. The interpretation of \citet{kovalev2024linear} combined with the GRAAL step in \cref{eq:extrapolation} suggests that one should choose $\ox_{k+1} = \alpha_k \hx_k + (1-\alpha_k)\ox_k$. However, this would require that the parameters $\alpha_k$ satisfy the following inequality:\footnote{More details on \cref{eq:alpha_eta} are available in the works of \citet{kovalev2024linear,kovalev2025sgd}.}
\begin{equation}\label{eq:alpha_eta}
  \eta_k/\alpha_k \leq \eta_{k-1}/\alpha_{k-1} + \eta_k.
\end{equation}
The best option is to choose $\alpha_k$ such that \cref{eq:alpha_eta} becomes an equality. However, this is impossible: computing $\eta_k$ requires an estimate of the local curvature, which requires the computation of the gradient $\nabla f_k(\hx_k)$ and its use in \cref{eq:Lambda}, which in turn requires knowing $\alpha_k$ in advance. Alternatively, one could follow the approach of \citet{li2025simple} and \citet{suh2025adaptive} used in AC-FGM and AdaNAG, respectively, and simply predefine $\alpha_k \propto 2/(k+2)$, just like in AGD. However, this requires additional restrictions on the stepsize $\eta_k$ and vastly limits the adaptation capabilities of AC-FGM, as we will discuss in \Cref{sec:L:comparison}.

\textbf{Solution: additional coupling step.} The key idea to resolve the issue above is to avoid the inequality in \cref{eq:alpha_eta} by defining $\ox_{k+1}$ differently, using an additional coupling step:
\begin{equation}
  \ox_{k+1} = \beta_k \tx_k + (1-\beta_k)\ox_k,\qquad
  \tx_k = \alpha_k \hx_k + (1-\alpha_k)\ox_k,
\end{equation}
where $\beta_k \in (0,1]$. Consequently, instead of requiring the inequality \cref{eq:alpha_eta}, we choose the parameters $\beta_k$ to satisfy the following relation:\footnote{More details on \cref{eq:alpha_beta_eta} are available in the proof of \Cref{thm:main} in \Cref{proof:thm:main}.}
\begin{equation}\label{eq:alpha_beta_eta}
  \eta_k/(\alpha_k \beta_k) = \eta_{k-1}/(\alpha_{k-1}\beta_{k-1}) + \eta_k = \cdots = H_k,
  \quad\text{where}\;\; H_k = \tsum_{i=0}^k \eta_i.
\end{equation}
Hence, we choose $\beta_k = \eta_k/(\alpha_kH_k)$ and avoid additional restrictions on the stepsize $\eta_k$. The only remaining question is how to ensure $\beta_k \leq 1$. The answer is that we choose $\alpha_k = \frac{(1+\gamma)\eta_{k-1}}{H_{k-1}+(1+\gamma)\eta_{k-1}}$. Indeed, in \Cref{lem:alpha_beta}, we prove that $\beta_k \in (0,1]$ by utilizing the inequality $\eta_{k}\leq (1+\gamma)\eta_{k-1}$, which is implied by our stepsize rule in \cref{eq:eta}. Moreover, our choice of $\alpha_k$ is implementable as it does not require knowledge of $\eta_k$, and, in contrast to AC-FGM and AdaNAG, it is adaptive because it is not based on any predefined sequence, but rather uses the adaptive stepsizes $\eta_{k-1}$ and their sum $H_{k-1}$.

\textbf{Adaptive stepsize.}
We use the following adaptive stepsize $\eta_k$ in our algorithm:
\begin{equation}\label{eq:eta}
  \eta_{k+1} = \min\brf*{(1+\gamma) \eta_k,\tfrac{\nu H_{k-1}\lambda_{k+1}}{\eta_{k-1}}},
\end{equation}
where $\lambda_{k+1} = \min\brf{\Lambda(\ox_{k+1};\tx_k),\Lambda(\ox_{k+1};\tx_{k+1})}$ is the local curvature estimator. This rule is primarily implied by the convergence analysis in the proof of \Cref{thm:main} in \Cref{proof:thm:main}. It can also be seen as a generalization of the stepsize rules in \cref{eq:GRAAL,eq:AdGD} for GRAAL and AdGD, respectively.

\begin{algorithm}[t]
  \caption{Accelerated GRAAL}
  \label{alg}
  \begin{algorithmic}[1]
    \State \textbf{input:} $x_0 \in \sX$, $\eta_0 > 0$, $K \in \brf{1,2,\dots}$
    \State \textbf{parameters:} $\theta,\gamma,\nu > 0$ satisfying \cref{eq:nu_gamma}
    \State $\alpha_0 = \beta_0 = 1$,\quad $H_0 = H_{-1} = \eta_{-1} = \eta_0$,\quad $\tx_0 = \ox_0 = x_0$
    \label{line:init}
    \For{$k=0,1,\dots,K-1$}
    \State $\alpha_{k+1} = \frac{(1+\gamma)\eta_k}{H_k + (1+\gamma)\eta_k}$
    \label{line:alpha}
    \State $x_{k+1} = x_k - \eta_k \nabla f(\tx_k)$\Comment{gradient step}
    \label{line:x}
    \State $\ox_{k+1} = \beta_k \tx_k + (1-\beta_k)\ox_k$\Comment{additional coupling step}
    \label{line:ox}
    \State $\hx_{k+1} = x_{k+1} + \theta(x_{k+1} - x_k)$\Comment{GRAAL extrapolation}
    \label{line:hx}
    \State $\tx_{k+1} = \alpha_{k+1}\hx_{k+1} + (1-\alpha_{k+1})\ox_{k+1}$\Comment{Nesterov acceleration/STM}
    \label{line:tx}
    \State $\lambda_{k+1} = \min\{\Lambda(\ox_{k+1};\tx_k),\Lambda(\ox_{k+1};\tx_{k+1})\}$\Comment{local curvature estimator}
    \label{line:lambda}
    \State $\eta_{k+1} = \min\brf*{(1+\gamma) \eta_k,\frac{\nu H_{k-1}\lambda_{k+1}}{\eta_{k-1}}}$,\;
    $H_{k+1} = H_k + \eta_{k+1}$
    \Comment{adaptive stepsize}
    \label{line:eta_H}
    \State $\beta_{k+1} = \frac{\eta_{k+1}}{\alpha_{k+1}H_{k+1}}$
    \label{line:beta}
    \EndFor
    \State \textbf{output:} $\ox_{K} \in \sX$
  \end{algorithmic}
\end{algorithm}

\subsection{Convergence Analysis}

We start the convergence analysis with the following two lemmas. In \Cref{lem:alpha_beta}, we show that $\beta_k \in (0,1]$ as discussed in the previous \Cref{sec:alg:dev}. In \Cref{lem:beta_f}, we use the additional coupling step from \cref{line:ox} and the convexity of the objective function $f(x)$ to obtain some useful inequalities.
\begin{lemma}<lem:alpha_beta>
  $\lambda_k,\eta_k, H_k > 0$, $\alpha_k \in (0,1)$, and $\beta_k \in (0,1]$ for all $k \in \brf{1,\dots,K}$.
\end{lemma}
\begin{lemma}<lem:beta_f>
  The following inequalities hold for all $k \in \brf{1,\dots,K-1}$:
  \begin{equation}
    f(\ox_k) - f(\tx_k) \leq \tfrac{1}{\beta_k}\brr{f(\ox_k) - f(\ox_{k+1})}
    ,\quad
    \bg_f(\ox_k;\tx_{k-1}) \leq \bg_f(\ox_{k-1};\tx_{k-1}).
  \end{equation}
\end{lemma}

Now, we obtain the main convergence result in \Cref{thm:main}. Note that \Cref{alg} requires the universal constant parameters $\theta,\gamma,\nu > 0$ to satisfy \cref{eq:nu_gamma}. However, it is easy to verify that such parameters exist. Using \Cref{thm:main}, we also obtain the upper bound on the functional suboptimality and distance to the solution of \cref{eq:main} in \Cref{cor:main}.
\begin{theorem}<thm:main>
  Let parameters $\theta,\gamma,\nu > 0$ satisfy the following relations:
  \begin{equation}\label{eq:nu_gamma}
    4\nu\theta(1+\gamma)^2 = \gamma,
    \quad
    1 + 2\gamma + \tfrac{2\gamma\theta^2}{(1+\theta)^2} \leq \tfrac{\theta}{(1+\theta)} + \tfrac{\theta^2}{(1+\theta)^2}.
  \end{equation}
  Then, the following inequality holds for all $x \in \sX$ and $k \in \brf{1,\dots,K-1}$:
  \begin{equation}
    \Psi_{k+1}(x) \leq \Psi_k(x)
    -\tfrac{\gamma\theta}{2}\eta_k^2\sqn{\nabla f(\tx_k)}
    -\tfrac{1}{4(1+\gamma)}\eta_k\bg_f(\ox_k;\tx_k),
  \end{equation}
  where $\Psi_k(x)$ is defined as follows:
  \begin{equation}\label{eq:Psi}
    \Psi_k(x) = \tfrac{1}{2}\sqn{x_k - x}
    +H_{k-1}\brr{f(\ox_k)-f(x)}
    +\tfrac{\theta\eta_k\eta_{k-1}}{\lambda_k}\bg_f(\ox_{k-1};\tx_{k-1})
    +\tfrac{\gamma\theta}{2}\sqn{x_k - x_{k-1}}.
  \end{equation}
\end{theorem}
\begin{corollary}<cor:main>
  The following inequality holds for all $x \in \sX$ and $K \geq 1$:
  \begin{equation}
    \tfrac{1}{2}\sqn{x_K - x}
    +H_{K-1}\brr{f(\ox_K)-f(x)}
    \leq
    \tfrac{1}{2}\sqn{x_0 - x}
    +\tfrac{(1+\gamma\theta)}{2}\eta_0^2\sqn{\nabla f(x_0)}.
    % -\eta_0\bg_f(x,x_0).
  \end{equation}
\end{corollary}
It is important to highlight that the results in \Cref{thm:main} and \Cref{cor:main} are general and do not require any additional assumptions on the objective function other than convexity and continuous differentiability. Consequently, they do not imply any non-asymptotic convergence results for \Cref{alg}. However, in \Cref{sec:L,sec:L01}, we will establish particular iteration complexities in the cases where the objective function is $L$-smooth and $(L_0,L_1)$-smooth, respectively.

\section[Convergence Analysis for L-Smooth Functions]{Convergence Analysis for $L$-Smooth Functions}\label{sec:L}

\subsection{Main Result}

In this section, we establish the iteration complexity of \Cref{alg} in the case where the objective function $f(x)$ is $L$-smooth for $L > 0$. That is, the gradient $\nabla f(x)$ is $L$-Lipschitz:
\begin{equation}
  \norm{\nabla f(x) - \nabla f(z)} \leq L \norm{x-z}
  \quad\text{for all}\;\; x,z\in\sX.
\end{equation}
We start with the following two lemmas. \Cref{lem:Lambda_L} provides a lower bound on the curvature estimates $\lambda_k$. This lemma is standard, and its proof is given, for instance, by \citet[Theorem~2.1.5]{nesterov2018lectures}. \Cref{lem:H} bounds the growth of the cumulative sum $H_k$ of the stepsizes $\eta_k$.
\begin{lemma}<lem:Lambda_L>
  $\lambda_k \geq 1/L$ for all $k \in \brf{1,\dots,K}$.
  % Let the function $f(x)$ have $L$-Lipschitz gradient. Then $\lambda_k \geq 1/L$ for all $k \in \brf{1,\dots,K}$.
\end{lemma}
\begin{lemma}<lem:H>
  $H_{k-1} \leq H_{k} \leq (2+\gamma)H_{k-1}$ for all $k \in \brf{0,\dots,K}$.
\end{lemma}

Now, we are ready to establish the lower bound on the cumulative sum $H_k$ of the stepsizes $\eta_k$ in the following \Cref{thm:L}. Using this bound, we establish the iteration complexity of \Cref{alg} for $L$-smooth functions in \Cref{cor:L}.
\begin{theorem}<thm:L>
  Let constants $c,m \in \R$ be defined as follows:
  \begin{equation}\label{eq:cm_L}
    c = \min\brf*{\tfrac{\sqrt{\nu}}{3(2+\gamma)},\tfrac{\sqrt{\nu}\gamma}{16\sqrt[4]{\gamma(1+\gamma)^5(2+\gamma)^3}}},\quad
    m = \ceil*{\max\brf*{2,\ln_{1+\gamma}\brs*{\tfrac{4c^2}{\gamma\eta_0 L}}}}.
  \end{equation}
  Then, the following inequality holds for all $k \in \brf{0,\dots,K}$:
  \begin{equation}\label{eq:H_L}
    \sqrt{H_k} \geq \tfrac{c}{\sqrt{L}}\cdot(k - m).
  \end{equation}
\end{theorem}

\begin{corollary}<cor:L>
  Let $\eta_0L \leq 1$.
  Then, to reach the precision $f(\ox_K) - f(x^*) \leq \epsilon$, the following number of iterations of \Cref{alg} is sufficient:
  \begin{equation}
    K = \cO\brr*{1 +  \sqrt{L\sqn{x_0 - x^*}/\epsilon} + \ln \brs*{\tfrac{1}{\eta_0L}}}.
  \end{equation}
\end{corollary}

Note that the complexity result in \Cref{cor:L} requires the initial stepsize $\eta_0$ to satisfy the inequality $\eta_0 \leq 1/L$. However, we can simply choose $\eta_0$ to be very small, say $10^{-10}$, as suggested by \citet{malitsky2020adaptive} for AdGD. This will only result in a small logarithmic additive factor in the iteration complexity as implied by \Cref{cor:L}.

\subsection{Comparison with AC-FGM and AdaNAG}\label{sec:L:comparison}

\textbf{AC-FGM.}
\citet[Corollary~1]{li2025simple} use the following adaptive stepsize rule for AC-FGM:
\begin{equation}\label{eq:AC-FGM}
  \eta_{k+1} = \min\brf*{\tfrac{k+1}{k}\eta_k, \tfrac{1}{8}k\lambda_{k+1}}.
\end{equation}
This rule implies that the stepsize growth is restricted by the inequality $\eta_{k+1} \leq (1+1/k)\eta_k$. As mentioned in \Cref{sec:alg:dev}, this restriction substantially limits the ability of AC-FGM to adapt to the curvature of the objective function. In particular, the result of \citet[Corollary~1]{li2025simple} implies the following iteration complexity, provided $\eta_0 \leq 0.4/L$:
\begin{equation}\label{eq:rate_AC-FGM}
  K = \cO\brr*{\sqrt{\max\brf{1,1/(\eta_0L)} \cdot L\sqn{x_0 - x^*}/\epsilon}}
  \quad\Rightarrow\quad
  f(\ox_K) - f(x^*) \leq \epsilon.
\end{equation}
This matches the optimal complexity in \cref{eq:rate_AGD}, if we choose $\eta_0 = 0.4/L$. However, if we choose the initial stepsize to be too small, i.e., $\eta_0 \ll 0.4/L$, the complexity result in \cref{eq:rate_AC-FGM} will be worse than the optimal one by a factor of $1/\sqrt{\eta_0L}$. In other words, the stepsize rule in \cref{eq:AC-FGM} cannot adapt to a ``bad'' choice of the initial stepsize $\eta_0$ due to stepsize growth restrictions. \citet{li2025simple} even had to use a line search at the first iteration of AC-FGM to find a ``good'' initial stepsize $\eta_0$ and achieve the optimal complexity in \cref{eq:rate_AGD}.

\textbf{AdaNAG.} \citet{suh2025adaptive} uses a stepsize rule in AdaNAG, which is not substantially different from \cref{eq:AC-FGM}. Consequently, they encounter issues similar to AC-FGM. The difference is that they do not use a line search at the first iteration, but rather estimate $\eta_0$ using Option~I in \cref{eq:options}, which implies the following iteration complexity:
\begin{equation}\label{eq:rate_AdaNAG}
  K = \cO\brr*{\max\brf{1,\eta_0L} \cdot \sqrt{L\sqn{x_0 - x^*}/\epsilon}}
  \quad\Rightarrow\quad
  f(\ox_K) - f(x^*) \leq \epsilon.
\end{equation}
This result may be significantly worse than the optimal one in \cref{eq:rate_AGD} if the initial stepsize estimate is too large, i.e., $\eta_0 L \gg 1$. Moreover, similar to AC-FGM, the growth of the stepsize in AdaNAG is also substantially restricted, which can limit its performance, for instance, under the more realistic $(L_0,L_1)$-smoothness assumption.

\textbf{\Cref{alg} vs AC-FGM and AdaNAG.} In contrast to AC-FGM and AdaNAG, our stepsize rule in \cref{eq:eta} allows the geometric growth of the stepsize $\eta_{k+1} \leq (1+\gamma)\eta_k$. Hence, \Cref{alg} can adapt even to a very small choice of the initial stepsize $\eta_0$ at the cost of a small logarithmic additive factor in the iteration complexity, as indicated by \Cref{cor:L}. In addition, as we will discuss in \Cref{sec:L01}, the geometric growth of the stepsize is crucial for adaptation under the $(L_0,L_1)$-smoothness assumption, where the local gradient Lipschitz constant may change at an exponential rate. It is also worth mentioning that \citet[Corollary~2]{li2025simple} and \citet[Theorem~6]{suh2025adaptive} tried to resolve the issues with the stepsize growth restrictions in AC-FGM and AdaNAG by using different stepsize rules. However, they could not properly justify the efficiency of these new stepsize rules and provably achieve geometric growth of the stepsize.

\section[Convergence Analysis for (L0,L1)-Smooth Functions]{Convergence Analysis for $(L_0,L_1)$-Smooth Functions}\label{sec:L01}

In this section, we establish the iteration complexity of \Cref{alg} in the case where the objective function $f(x)$ is $(L_0,L_1)$-smooth for $L_0,L_1 > 0$. That is, the objective function $f(x)$ is twice continuously differentiable and the following inequality holds:
\begin{equation}
  \norm{\nabla^2 f(x)} \leq L_0 + L_1\norm{\nabla f(x)}
  \;\;\text{for all}\;\; x \in \sX.
\end{equation}
This assumption was proposed by \citet{zhang2019gradient} and is primarily motivated by experiments suggesting that the norm of the Hessian correlates with the gradient norm of the objective functions in deep neural networks. Note that the requirement for twice continuous differentiability may be relaxed by using the following equivalent condition for continuously differentiable objective functions \citep[Lemma~2.5]{vankov2024optimizing}:
\begin{equation}\label{eq:L01}
  \norm{\nabla f(x) - \nabla f(z)} \leq \brr{L_0 + L_1\norm{\nabla f(x)}} \cdot \tfrac{1}{L_1}\brr{\exp\brr{L_1\norm{x-z}} - 1}.
\end{equation}
It is also important to highlight that $(L_0,L_1)$-smoothness implies $L$-smoothness with $L=L_0$. The reverse is obviously not true: $(L_0,L_1)$-smoothness is much more general and allows the exponential growth of the objective function and local gradient Lipschitz constant (\citealp[Lemma~2.1]{gorbunov2024methods}; \citealp[Lemma~2.5]{vankov2024optimizing}).

\subsection{Main Result}\label{sec:L01:main}

We start the convergence analysis of \Cref{alg} with \Cref{lem:D}, which refines the previously obtained results in \Cref{thm:main} and \Cref{cor:main}. Furthermore, in \Cref{lem:lambda_min,lem:lambda_k}, we establish lower bounds on the estimate $\lambda_k$ of the local inverse gradient Lipschitz constant.
\begin{lemma}<lem:D>
  The following inequalities hold for all $K \geq 1$:
  \begin{equation}
    \norm{\tx_{K-1} - x^*}\leq \cD,\quad
    \norm{\ox_{K-1} - x^*} \leq \cD,\quad
    \tsum_{i=1}^K\brr*{\eta_i^2\sqn{\nabla f(\tx_i)} + \eta_i\bg_f(\ox_i;\tx_i)}
    \leq
    \cD^2,
  \end{equation}
  where $\cD \geq 0$ is defined as follows:
  \begin{equation}\label{eq:D}
    \cD^2 = \max\brf*{\tfrac{1}{\gamma\theta},2(1+\gamma), (1+2\theta)^2}
    \brr*{
      \sqn{x_0 - x^*}
    +(1+\gamma\theta)\eta_0^2\sqn{\nabla f(x_0)}}.
  \end{equation}
\end{lemma}

\begin{lemma}<lem:lambda_min>
  $\lambda_k \geq \lmin$ for all $k \in \brf{1,\dots,K}$, where $\lmin > 0$ is defined as follows:
  \begin{equation}
    \lmin = \tfrac{1}{L_0}\exp(-3L_1\cD).
  \end{equation}
\end{lemma}

\begin{lemma}<lem:lambda_k>
  The following inequality holds for all $k \in \brf{1,\dots,K}$:
  \begin{equation}
    \resizebox{0.93\textwidth}{!}{
      $\displaystyle\lambda_k \geq \frac{4}{9\max\brf*{2L_0, 2L_1\norm{\nabla f(\tx_k)},2L_1\norm{\nabla f(\tx_{k-1})},9L_1^2\bg_f(\ox_k,\tx_k),9L_1^2\bg_f(\ox_{k-1},\tx_{k-1})}}.$
    }
  \end{equation}
\end{lemma}

Next, we define the sets of indices $\cT_1(k),\cT_2(k),\cT_3(k),\cT_4(k)$ as follows:
\begin{equation}
  \begin{aligned}\label{eq:Tau}
    \cT_1(k) &= \brf*{1\leq i \leq k :\eta_i = \tfrac{\nu H_{i-2}\lambda_i}{\eta_{i-2}}\text{\;and\;}\lambda_i \geq \tfrac{2}{9L_0}},\\
    \cT_2(k) &= \brf*{1\leq i \leq k :\eta_i = \tfrac{\nu H_{i-2}\lambda_i}{\eta_{i-2}}\text{\;and\;}\lambda_i < \tfrac{2}{9L_0}},\\
    \cT_3(k) &= \brf*{1\leq i \leq k :\eta_i < \tfrac{\nu H_{i-2}\lambda_i}{\eta_{i-2}}\text{\;and\;} \brr*{l(i) \in \cT_1(k) \text{\;or\;} i - l(i) > m}},\\
    \cT_4(k) &= \brf*{1\leq i \leq k :\eta_i < \tfrac{\nu H_{i-2}\lambda_i}{\eta_{i-2}}\text{\;and\;} l(i) \in \cT_2(k)\cup\brf{0} \text{\;and\;} i - l(i) \leq m},\\
  \end{aligned}
\end{equation}
where $m$ is a positive integer, and the integer function $l(k)$ is defined as follows:
\begin{equation}\label{eq:l}
  l(k) = \max \brf*{i : i \in \cT_1(k)\cup \cT_2(k) \cup \brf{0}}.
\end{equation}
It is not hard to verify that these sets of indices are pairwise disjoint and that $\cup_{j=1}^4 \cT_j(k) = \brf{1,\dots,k}$. Moreover, the sizes of the sets $\cT_2(k)$ and $\cT_4(k)$ are bounded as shown in the following \Cref{lem:Tau}.

\begin{lemma}<lem:Tau>
  The following inequalities hold for all $k \in \brf{1,\dots,K}$:
  \begin{equation}
    \abs{\cT_4(k)} \leq m + m \abs{\cT_2(k)}
    ,\qquad
    \abs{\cT_2(k)} \leq 1 + \tfrac{81(1+\gamma)^2}{2\min\brf{\nu,\nu^2}}\cdot L_1^2\cD^2.
  \end{equation}
\end{lemma}

Now, we establish the key lower bound on the cumulative sum $H_k$ of the stepsizes $\eta_k$ in the following \Cref{thm:L01}. Using this bound, we establish the iteration complexity of \Cref{alg} for $(L_0,L_1)$-smooth functions in \Cref{cor:L01}.

\begin{theorem}<thm:L01>
  Let $c > 0$ and $m > 0$  be defined as follows:
  \begin{equation}\label{eq:cm_L01}
    \resizebox{0.93\textwidth}{!}{$\displaystyle
      c = \min\brf*{ \tfrac{\sqrt{2\nu}}{9(2+\gamma)},\tfrac{\sqrt{\nu}\gamma}{24\sqrt[4]{4\gamma(1+\gamma)^5(2+\gamma)^3}}},
      \quad
      m = \ceil*{\max\brf*{1,\ln_{1+\gamma}\brs*{\tfrac{4c^2}{\gamma \eta_0 L_0}},4\ln_{1+\gamma}\brs*{\tfrac{2}{9L_0\lmin}}}}.
    $}
  \end{equation}
  Then, the following inequality holds for all $k \in \brf{0,\dots,K}$.
  \begin{equation}\label{eq:H_L01}
    \sqrt{H_k} \geq \tfrac{c}{\sqrt{L_0}} \cdot \brr*{k - \abs{\cT_2(k)} - \abs{\cT_4(k)} - 1}.
  \end{equation}
\end{theorem}

% \begin{corollary}<cor:L01>
%   To reach the precision $f(\ox_K) - f(x^*) \leq \epsilon$, it is sufficient to perform the following number of iterations of \Cref{alg}:
%   \begin{equation}\label{eq:K_L01}
%     K = \cO\brr*{1 + \sqrt{L_0\cD^2/\epsilon} + L_1^3 \cD^3 + \brr*{1 + L_1^2\cD^2}\ln \brs*{\tfrac{1}{\eta_0L_0}}}.
%   \end{equation}
%   Moreover, choosing the initial stepsize $\eta_0 \leq \norm{x_0 - x^*} / \norm{\nabla f(x_0)}$ implies $\cD = \cO(\norm{x_0 - x^*})$.
% \end{corollary}
\begin{corollary}<cor:L01>
  Let $\eta_0L_0\exp(L_1\norm{x_0-x^*}) \leq 1$.
  Then, $\cD = \cO(\norm{x_0 - x^*})$, and to reach the precision $f(\ox_K) - f(x^*) \leq \epsilon$, the following number of iterations of \Cref{alg} is sufficient:
  \begin{equation}\label{eq:K_L01}
    K = \cO\brr*{1 + \sqrt{L_0\cD^2/\epsilon} + L_1^3 \cD^3 + \brr*{1 + L_1^2\cD^2}\ln \brs*{\tfrac{1}{\eta_0L_0}}}.
  \end{equation}
\end{corollary}

Similar to \Cref{cor:L} for the $L$-smooth case, \Cref{cor:L01} requires the initial stepsize $\eta_0$ to satisfy the inequality $\eta_0L_0\exp(L_1\norm{x_0-x^*})\leq 1$. We can ensure this inequality without any line search or hyperparameter tuning, simply by choosing a very small initial stepsize $\eta_0$. Choosing the initial stepsize $\eta_0$ too small will only result in an additive constant factor $\brr*{1 + L_1^2\cD^2}\ln \big[\frac{1}{\eta_0L_0}\big]$, which does not depend on the precision $\epsilon$ and has a logarithmic dependence on $\eta_0$.

\subsection{Comparison with Existing Results}\label{sec:L01:comparison}

\begin{table}[t]
  \caption{Comparison of the iteration complexities for solving \cref{eq:main} under the convexity and $(L_0,L_1)$-smoothness; universal constants are omitted; $\cD = \norm{x_0 - x^*}$; the initial functional gap is bounded as $f(x_0) - f(x^*) \leq \cO\brr{L_0\cD^2\exp(L_1\cD)}$, where necessary; optimality is considered up to additive constants.}
  \label{tab}
  \centering
  \begin{NiceTabular}{|c|l|c|c|}
    \CodeBefore
    \rowcolor{cyan!10}{6}
    \Body
    \toprule
    \bf Reference & \Block[c]{1-1}{\bf Iteration Complexity} & \bf Optimal & \bf Adaptive
    \\\midrule
    \citet{li2023convex} & ${\color{red}\sqrt{\frac{L_0\cD^2}{\epsilon}}} \times \brr*{1 + L_1^2\cD^2}\exp(\cO(1)L_1\cD)$ & \xmark & \xmark
    \\\midrule
    \citet{gorbunov2024methods} & ${\color{red}\sqrt{\frac{L_0\cD^2}{\epsilon}}} \times \sqrt{1 + L_1\cD\exp(L_1\cD)}$  & \xmark & \xmark
    \\\midrule
    \citet{vankov2024optimizing} & ${\color{red}\sqrt{\frac{L_0\cD^2}{\epsilon}}} + (L_1\cD)^{5/3}$ &\cmark & \xmark
    \\\midrule
    \citet{tyurin2025near} & ${\color{red}\sqrt{\frac{L_0\cD^2}{\epsilon}}} + (L_1\cD)^{2}$ &\cmark & \xmark
    \\\midrule
    \bf\Cref{cor:L01} & ${\color{red}\sqrt{\frac{L_0\cD^2}{\epsilon}}} + (L_1\cD)^{3}$ &\cmark & \cmark
    \\\bottomrule
  \end{NiceTabular}
\end{table}

\textbf{Accelerated methods for $(L_0,L_1)$-smooth functions.} As mentioned in \Cref{sec:related}, there are several existing accelerated algorithms with theoretical guarantees for minimizing convex $(L_0,L_1)$-smooth functions. To compare these results with \Cref{alg}, we use a particular choice of the initial stepsize $\eta_0 = \frac{1}{L_0}\exp(-L_1\norm{x_0-x^*})$ in \Cref{cor:L01}. The comparison is summarized in \Cref{tab}. The algorithms of \citet{li2023convex,gorbunov2024methods} are neither adaptive nor optimal. The algorithms of \citet{vankov2024optimizing,tyurin2025near} are near-optimal as they match the complexity in \cref{eq:rate_AGD} up to additive constants, just like \Cref{alg}. The result of \citet{vankov2024optimizing} has a slightly better additive constant $(L_1\cD)^{5/3}$. However, neither of the algorithms of \citet{vankov2024optimizing,tyurin2025near} is adaptive: the algorithm of \citet{vankov2024optimizing} requires solving a one-dimensional auxiliary optimization subproblem at each iteration, and the algorithm of \citet{tyurin2025near} requires tuning several parameters. In contrast, \Cref{alg} does not require any hyperparameter tuning or line search to achieve near-optimal complexity, as discussed in the previous \Cref{sec:L01:main}.

\textbf{AdGD for $(L_0,L_1)$-smooth functions.}
\citet{gorbunov2024methods} established the iteration complexity $\cO\brr*{L_0\cD^2/\epsilon + (L_1\cD)^6}$ for the AdGD algorithm under the $(L_0,L_1)$-smoothness assumption.\footnote{Note that \citet{gorbunov2024methods} incorrectly reported their result as $\tilde\cO\brr*{L_0\cD^2/\epsilon + (L_1\cD)^4}$.} This result is unsurprisingly worse than ours in \Cref{cor:L01} due to the lack of acceleration. In addition, \citet{gorbunov2024methods} did not prove that the constant $\cD$ is bounded as $\cD = \cO\brr{\norm{x_0 - x^*}}$. In fact, $\cD$ also contains the initial objective function gap $f(x_0) - f(x^*)$, and hence, it may have an exponential dependency on the initial distance $\norm{x_0 - x^*}$.

\textbf{AC-FGM and AdaNAG.}
As previously discussed in \Cref{sec:L:comparison}, we allow the geometric growth of the adaptive stepsize in \Cref{alg}, which is crucial for obtaining the near-optimal complexity result in \Cref{cor:L01}. Indeed, the estimates of the local curvature $\lambda_k$ can scale exponentially in the worst case according to \Cref{lem:lambda_min}, but can grow up to $\cO(1/L_0)$ when the algorithm reaches a certain region near the solution $x^*$. Hence, the growth of the stepsize at a geometric rate or faster is necessary to avoid exponential factors in the iteration complexity. In contrast, \citet{li2025simple,suh2025adaptive} do not provide any convergence guarantees for $(L_0,L_1)$-smoothness for AC-FGM and AdaNAG. In addition, as previously discussed, the rate of stepsize growth in AC-FGM and AdaNAG is far below geometric. Hence, we conjecture that it is not possible to reach near-optimal complexity with these algorithms.

\newpage

\cite*{}

\bibliographystyle{apalike}
\bibliography{references}

\newpage

\appendix

\part*{Appendix}

\section[Proofs for \crtCref{sec:alg}]{Proofs for \Cref{sec:alg}}

\proofsubsection{lem:alpha_beta}

First, we can show that $\lambda_k > 0$ for all $k \in \brf{1,\dots,K}$. Otherwise, if $\lambda_k = 0$, there exist $x,z \in \sX$ such that $\Lambda(x,z) = 0$ due to the definition of $\lambda_k$ on \cref{line:lambda}. Hence, $\bg_f(x;z) = 0$ according to \cref{eq:Lambda}. From the convexity of the function $f(x)$, it follows that $\bg_f(x;z) \geq 0$. Hence, $x$ minimizes the function $\bg_f(\cdot;z)$, and by the first-order optimality conditions, $\nabla_x \bg_f(x;z) = \nabla f(x) - \nabla f(z) = 0$. This implies $\Lambda(x;z) = +\infty$ according to \cref{eq:Lambda}, which contradicts the possibility of $\lambda_k = 0$.

Next, using the fact that $\lambda_k > 0$ and \cref{line:eta_H,line:init}, it is easy to show that $\eta_k, H_k > 0$ for all $k \in \brf{-1,\dots,K}$. Hence, the inclusions $\alpha_k \in (0,1)$ and $\beta_k \in (0,+\infty)$ are obvious due to their definitions on \cref{line:alpha,line:beta}. Finally, we can upper-bound $\beta_{k+1}$ as follows:
\begin{align*}
  \beta_{k+1}
  &\at{uses \cref{line:beta}}=
  \frac{\eta_{k+1}}{\alpha_{k+1}H_{k+1}}
  \at{uses \cref{line:alpha,line:eta_H}}=
  \frac{\eta_{k+1}}{H_k + \eta_{k+1}}\cdot\frac{H_k + (1+\gamma)\eta_k}{(1+\gamma)\eta_k}
  \at{uses the inequality $\eta_{k+1}\leq (1+\gamma)\eta_k$ implied by \cref{line:eta_H} and the monotonicity of the function $t \mapsto t/(1+t)$}\leq
  \frac{(1+\gamma)\eta_k}{H_k + (1+\gamma)\eta_k}\cdot\frac{H_k + (1+\gamma)\eta_k}{(1+\gamma)\eta_k}
  =1,
\end{align*}
where \annotate.\qed

\proofsubsection{lem:beta_f}

We can upper-bound $f(\ox_{k+1})$ as follows:
\begin{align*}
  f(\ox_{k+1})
  &\at{uses \cref{line:ox}}=
  f(\beta_k \tx_k + (1-\beta_k)\ox_k)
  \at{uses \Cref{lem:alpha_beta} and the convexity of $f(x)$}\leq
  \beta_k f(\tx_k) + (1-\beta_k)f(\ox_k),
\end{align*}
where\annotate. After rearranging, we obtain the first desired inequality.
Furthermore, we can upper-bound $\bg_f(\ox_k,\tx_{k-1})$ as follows:
\begin{align*}
  \bg_f(\ox_k;\tx_{k-1})
  &\at{uses \cref{line:ox}}=
  \bg_f(\beta_{k-1}\tx_{k-1} + (1-\beta_{k-1})\ox_{k-1};\tx_{k-1})
  \\&\at{use uses \Cref{lem:alpha_beta} and the convexity of $\bg_f(\cdot;\tx_{k-1})$, which is implied by the convexity of $f(x)$}\leq
  (1-\beta_{k-1})\bg_f(\ox_{k-1};\tx_{k-1}) + \beta_{k-1}\bg_f(\tx_{k-1};\tx_{k-1})
  \\&\at{use uses \Cref{lem:alpha_beta} and the convexity of $\bg_f(\cdot;\tx_{k-1})$, which is implied by the convexity of $f(x)$}\leq
  \bg_f(\ox_{k-1};\tx_{k-1})
\end{align*}
where \annotate. This proves the second desired inequality.\qed

\proofsubsection{thm:main}

For all $k \in \brf{1,\dots,K-1}$, we obtain the following:
\begin{align*}
  \tfrac{1}{2}\sqn{x_{k+1} - x}
  &=
  \tfrac{1}{2}\sqn{x_k - x}
  -\tfrac{1}{2}\sqn{x_{k+1} - x_k}
  +\<x_{k+1} - x_k,x_{k+1} - x>
  \\&=
  \tfrac{1}{2}\sqn{x_k - x}
  -\tfrac{1}{2}\sqn{x_{k+1} - x_k}
  +\<x_{k+1}-x_k,x_{k+1} - \hx_{k+1}>
  \\&
  +\<x_{k+1}-x_k,\hx_{k+1} - \hx_k + \hx_k - x>
  \\&\at{uses \cref{line:hx}}=
  \tfrac{1}{2}\sqn{x_k - x}
  -\brr*{\tfrac{1}{2} + \theta}\sqn{x_{k+1} - x_k}
  +\<x_{k+1}-x_k,\hx_{k+1} - \hx_k + \hx_k - x>
  \\&\at{uses \cref{line:x}}=
  \tfrac{1}{2}\sqn{x_k - x}
  -\brr*{\tfrac{1}{2} + \theta}\sqn{x_{k+1} - x_k}
  -\eta_k\<\nabla f(\tx_k),\hx_{k+1} - \hx_k>
  \\&
  +\eta_k\<\nabla f(\tx_k),x - \hx_k>
  \\&\at{uses \cref{line:tx}}=
  \tfrac{1}{2}\sqn{x_k - x}
  -\brr*{\tfrac{1}{2} + \theta}\sqn{x_{k+1} - x_k}
  -\eta_k\<\nabla f(\tx_k),\hx_{k+1} - \hx_k>
  \\&
  +\eta_k\<\nabla f(\tx_k),x - \tx_k>
  +\tfrac{(1-\alpha_k)\eta_k}{\alpha_k}\<\nabla f(\tx_k),\ox_k - \tx_k>
  \\&\at{uses the convexity of function $f(x)$}\leq
  \tfrac{1}{2}\sqn{x_k - x}
  -\brr*{\tfrac{1}{2} + \theta}\sqn{x_{k+1} - x_k}
  -\eta_k\<\nabla f(\tx_k),\hx_{k+1} - \hx_k>
  \\&
  +\eta_k\brr{f(x) - f(\tx_k)}
  +\tfrac{(1-\alpha_k)\eta_k}{\alpha_k}\brr{f(\ox_k) - f(\tx_k) - \bg_f(\ox_k;\tx_k)}
  \\&=
  \tfrac{1}{2}\sqn{x_k - x}
  -\brr*{\tfrac{1}{2} + \theta}\sqn{x_{k+1} - x_k}
  -\eta_k\<\nabla f(\tx_k),\hx_{k+1} - \hx_k>
  \\&
  -\eta_k \brr{f(\ox_k) - f(x)}
  +\tfrac{\eta_k}{\alpha_k}\brr{f(\ox_k)-f(\tx_k)}
  -\tfrac{(1-\alpha_k)\eta_k}{\alpha_k}\bg_f(\ox_k;\tx_k)
  \\&\at{uses \Cref{lem:beta_f}}\leq
  \tfrac{1}{2}\sqn{x_k - x}
  -\brr*{\tfrac{1}{2} + \theta}\sqn{x_{k+1} - x_k}
  -\eta_k\<\nabla f(\tx_k),\hx_{k+1} - \hx_k>
  \\&
  -\eta_k \brr{f(\ox_k) - f(x)}
  +\tfrac{\eta_k}{\alpha_k\beta_k}\brr{f(\ox_k)-f(\ox_{k+1})}
  -\tfrac{(1-\alpha_k)\eta_k}{\alpha_k}\bg_f(\ox_k;\tx_k)
  \\&\at{uses \cref{line:beta}}=
  \tfrac{1}{2}\sqn{x_k - x}
  -\brr*{\tfrac{1}{2} + \theta}\sqn{x_{k+1} - x_k}
  -\eta_k\<\nabla f(\tx_k),\hx_{k+1} - \hx_k>
  \\&
  -\eta_k \brr{f(\ox_k) - f(x)}
  +H_k\brr{f(\ox_k)-f(\ox_{k+1})}
  -\tfrac{(1-\alpha_k)\eta_k}{\alpha_k}\bg_f(\ox_k;\tx_k)
  \\&\at{uses \cref{line:eta_H}}=
  \tfrac{1}{2}\sqn{x_k - x}
  -\brr*{\tfrac{1}{2} + \theta}\sqn{x_{k+1} - x_k}
  -\eta_k\<\nabla f(\tx_k),\hx_{k+1} - \hx_k>
  \\&
  +H_{k-1}\brr{f(\ox_k)-f(x)}
  -H_k\brr{f(\ox_{k+1}) - f(x)}
  -\tfrac{(1-\alpha_k)\eta_k}{\alpha_k}\bg_f(\ox_k;\tx_k),
\end{align*}
where \annotate.
Furthermore, we get the following:
\begin{align*}
  \mi{1}\tfrac{1}{2}\sqn{x_{k+1} - x}
  \\&\leq
  \tfrac{1}{2}\sqn{x_k - x}
  -\brr*{\tfrac{1}{2} + \theta}\sqn{x_{k+1} - x_k}
  +H_{k-1}\brr{f(\ox_k)-f(x)}
  -H_k\brr{f(\ox_{k+1}) - f(x)}
  \\&
  -\tfrac{(1-\alpha_k)\eta_k}{\alpha_k}\bg_f(\ox_k;\tx_k)
  -\eta_k\<\nabla f(\tx_k) - \nabla f(\tx_{k-1}),\hx_{k+1} - \hx_k>
  -\eta_k\<\nabla f(\tx_{k-1}),\hx_{k+1} - \hx_k>
  \\&\at{uses \cref{line:x}}=
  \tfrac{1}{2}\sqn{x_k - x}
  -\brr*{\tfrac{1}{2} + \theta}\sqn{x_{k+1} - x_k}
  +H_{k-1}\brr{f(\ox_k)-f(x)}
  -H_k\brr{f(\ox_{k+1}) - f(x)}
  \\&
  -\tfrac{(1-\alpha_k)\eta_k}{\alpha_k}\bg_f(\ox_k;\tx_k)
  -\eta_k\<\nabla f(\tx_k) - \nabla f(\tx_{k-1}),\hx_{k+1} - \hx_k>
  \\&
  +\tfrac{\eta_k}{\eta_{k-1}}\<x_k - x_{k-1},\hx_{k+1} - \hx_k>
  \\&\at{use \cref{line:hx}}=
  \tfrac{1}{2}\sqn{x_k - x}
  -\brr*{\tfrac{1}{2} + \theta}\sqn{x_{k+1} - x_k}
  +H_{k-1}\brr{f(\ox_k)-f(x)}
  -H_k\brr{f(\ox_{k+1}) - f(x)}
  \\&
  -\tfrac{(1-\alpha_k)\eta_k}{\alpha_k}\bg_f(\ox_k;\tx_k)
  -\eta_k\<\nabla f(\tx_k) - \nabla f(\tx_{k-1}),\hx_{k+1} - \hx_k>
  \\&
  +\tfrac{\eta_k}{\theta\eta_{k-1}}\<\hx_k - x_k,\hx_{k+1} - \hx_k>
  \\&=
  \tfrac{1}{2}\sqn{x_k - x}
  -\brr*{\tfrac{1}{2} + \theta}\sqn{x_{k+1} - x_k}
  +H_{k-1}\brr{f(\ox_k)-f(x)}
  -H_k\brr{f(\ox_{k+1}) - f(x)}
  \\&
  -\tfrac{(1-\alpha_k)\eta_k}{\alpha_k}\bg_f(\ox_k;\tx_k)
  -\eta_k\<\nabla f(\tx_k) - \nabla f(\tx_{k-1}),\hx_{k+1} - \hx_k>
  \\&
  +\tfrac{\eta_k}{2\theta\eta_{k-1}}\brr{\sqn{\hx_{k+1} - x_k} - \sqn{\hx_k - x_k} - \sqn{\hx_{k+1} - \hx_k}}
  \\&\at{use \cref{line:hx}}=
  \tfrac{1}{2}\sqn{x_k - x}
  -\brr*{\tfrac{1}{2} + \theta}\sqn{x_{k+1} - x_k}
  +H_{k-1}\brr{f(\ox_k)-f(x)}
  -H_k\brr{f(\ox_{k+1}) - f(x)}
  \\&
  -\tfrac{(1-\alpha_k)\eta_k}{\alpha_k}\bg_f(\ox_k;\tx_k)
  -\eta_k\<\nabla f(\tx_k) - \nabla f(\tx_{k-1}),\hx_{k+1} - \hx_k>
  \\&
  +\tfrac{\eta_k}{2\theta\eta_{k-1}}\brr{(1+\theta)^2\sqn{x_{k+1} - x_k} - \theta^2\sqn{x_k - x_{k-1}} - \sqn{\hx_{k+1} - \hx_k}}
  \\&\at{uses the Cauchy-Schwarz inequality}\leq
  \tfrac{1}{2}\sqn{x_k - x}
  +\tfrac{(1+\theta)^2}{2\theta}\brr*{\tfrac{\eta_k}{\eta_{k-1}} - \tfrac{\theta + 2\theta^2}{(1+\theta)^2}}\sqn{x_{k+1} - x_k}
  -\tfrac{\theta\eta_k}{2\eta_{k-1}}\sqn{x_k - x_{k-1}}
  \\&
  -\tfrac{\eta_k}{2\theta\eta_{k-1}}\sqn{\hx_{k+1} - \hx_k}
  +\eta_k\norm{\nabla f(\tx_k) - \nabla f(\tx_{k-1})}\norm{\hx_{k+1} - \hx_k}
  \\&
  -\tfrac{(1-\alpha_k)\eta_k}{\alpha_k}\bg_f(\ox_k;\tx_k)
  +H_{k-1}\brr{f(\ox_k)-f(x)}
  -H_k\brr{f(\ox_{k+1}) - f(x)},
\end{align*}
where \annotate.
Next, we get
\begin{align*}
  \mi{1}\tfrac{1}{2}\sqn{x_{k+1} - x}
  \\&\leq
  \tfrac{1}{2}\sqn{x_k - x}
  +\tfrac{(1+\theta)^2}{2\theta}\brr*{\tfrac{\eta_k}{\eta_{k-1}} - \tfrac{\theta + 2\theta^2}{(1+\theta)^2}}\sqn{x_{k+1} - x_k}
  -\tfrac{\theta\eta_k}{2\eta_{k-1}}\sqn{x_k - x_{k-1}}
  \\&
  -\tfrac{(1-\alpha_k)\eta_k}{\alpha_k}\bg_f(\ox_k;\tx_k)
  +H_{k-1}\brr{f(\ox_k)-f(x)}
  -H_k\brr{f(\ox_{k+1}) - f(x)}
  \\&
  -\tfrac{\eta_k}{2\theta\eta_{k-1}}\sqn{\hx_{k+1} - \hx_k}
  +\eta_k\norm{\nabla f(\tx_k) - \nabla f(\tx_{k-1})}\norm{\hx_{k+1} - \hx_k}
  \\&\at{uses the triangle inequality}\leq
  \tfrac{1}{2}\sqn{x_k - x}
  +H_{k-1}\brr{f(\ox_k)-f(x)}
  -H_k\brr{f(\ox_{k+1}) - f(x)}
  \\&
  +\tfrac{(1+\theta)^2}{2\theta}\brr*{\tfrac{\eta_k}{\eta_{k-1}} - \tfrac{\theta + 2\theta^2}{(1+\theta)^2}}\sqn{x_{k+1} - x_k}
  -\tfrac{\theta\eta_k}{2\eta_{k-1}}\sqn{x_k - x_{k-1}}
  \\&
  -\tfrac{\eta_k}{2\theta\eta_{k-1}}\sqn{\hx_{k+1} - \hx_k}
  -\tfrac{(1-\alpha_k)\eta_k}{\alpha_k}\bg_f(\ox_k;\tx_k)
  \\&
  +\eta_k\norm{\nabla f(\tx_k) - \nabla f(\ox_k)}\norm{\hx_{k+1} - \hx_k}
  +\eta_k\norm{\nabla f(\ox_k) - \nabla f(\tx_{k-1})}\norm{\hx_{k+1} - \hx_k}
  \\&\at{uses \cref{line:lambda} and \cref{eq:Lambda}}\leq
  \tfrac{1}{2}\sqn{x_k - x}
  +H_{k-1}\brr{f(\ox_k)-f(x)}
  -H_k\brr{f(\ox_{k+1}) - f(x)}
  \\&
  +\tfrac{(1+\theta)^2}{2\theta}\brr*{\tfrac{\eta_k}{\eta_{k-1}} - \tfrac{\theta + 2\theta^2}{(1+\theta)^2}}\sqn{x_{k+1} - x_k}
  -\tfrac{\theta\eta_k}{2\eta_{k-1}}\sqn{x_k - x_{k-1}}
  \\&
  -\tfrac{\eta_k}{2\theta\eta_{k-1}}\sqn{\hx_{k+1} - \hx_k}
  -\tfrac{(1-\alpha_k)\eta_k}{\alpha_k}\bg_f(\ox_k;\tx_k)
  \\&
  +\eta_k\sqrt{\tfrac{2}{\lambda_k}\bg_f(\ox_k;\tx_k)}\norm{\hx_{k+1} - \hx_k}
  +\eta_k\sqrt{\tfrac{2}{\lambda_k}\bg_f(\ox_k;\tx_{k-1})}\norm{\hx_{k+1} - \hx_k}
  \\&\at{uses \Cref{lem:beta_f}}\leq
  \tfrac{1}{2}\sqn{x_k - x}
  +H_{k-1}\brr{f(\ox_k)-f(x)}
  -H_k\brr{f(\ox_{k+1}) - f(x)}
  \\&
  +\tfrac{(1+\theta)^2}{2\theta}\brr*{\tfrac{\eta_k}{\eta_{k-1}} - \tfrac{\theta + 2\theta^2}{(1+\theta)^2}}\sqn{x_{k+1} - x_k}
  -\tfrac{\theta\eta_k}{2\eta_{k-1}}\sqn{x_k - x_{k-1}}
  \\&
  -\tfrac{\eta_k}{2\theta\eta_{k-1}}\sqn{\hx_{k+1} - \hx_k}
  -\tfrac{(1-\alpha_k)\eta_k}{\alpha_k}\bg_f(\ox_k;\tx_k)
  \\&
  +\eta_k\sqrt{\tfrac{2}{\lambda_k}\bg_f(\ox_k;\tx_k)}\norm{\hx_{k+1} - \hx_k}
  +\eta_k\sqrt{\tfrac{2}{\lambda_k}\bg_f(\ox_{k-1};\tx_{k-1})}\norm{\hx_{k+1} - \hx_k}
  \\&\at{uses Young's inequality}\leq
  \tfrac{1}{2}\sqn{x_k - x}
  +H_{k-1}\brr{f(\ox_k)-f(x)}
  -H_k\brr{f(\ox_{k+1}) - f(x)}
  \\&
  +\tfrac{(1+\theta)^2}{2\theta}\brr*{\tfrac{\eta_k}{\eta_{k-1}} - \tfrac{\theta + 2\theta^2}{(1+\theta)^2}}\sqn{x_{k+1} - x_k}
  -\tfrac{\theta\eta_k}{2\eta_{k-1}}\sqn{x_k - x_{k-1}}
  \\&
  -\tfrac{\eta_k}{2\theta\eta_{k-1}}\sqn{\hx_{k+1} - \hx_k}
  -\tfrac{(1-\alpha_k)\eta_k}{\alpha_k}\bg_f(\ox_k;\tx_k)
  +\tfrac{(1-\alpha_k)\eta_k}{2\alpha_k}\bg_f(\ox_k;\tx_k)
  \\&
  +\tfrac{\alpha_k\eta_k}{(1-\alpha_k)\lambda_k}\sqn{\hx_{k+1} - \hx_k}
  +\tfrac{\theta\eta_k\eta_{k-1}}{\lambda_k}\bg_f(\ox_{k-1};\tx_{k-1})
  +\tfrac{\eta_k}{2\theta\eta_{k-1}}\sqn{\hx_{k+1} - \hx_k}
  \\&=
  \tfrac{1}{2}\sqn{x_k - x}
  +H_{k-1}\brr{f(\ox_k)-f(x)}
  -H_k\brr{f(\ox_{k+1}) - f(x)}
  -\tfrac{(1-\alpha_k)\eta_k}{4\alpha_k}\bg_f(\ox_k;\tx_k)
  \\&
  +\tfrac{(1+\theta)^2}{2\theta}\brr*{\tfrac{\eta_k}{\eta_{k-1}} - \tfrac{\theta + 2\theta^2}{(1+\theta)^2}}\sqn{x_{k+1} - x_k}
  -\tfrac{\theta\eta_k}{2\eta_{k-1}}\sqn{x_k - x_{k-1}}
  \\&
  +\tfrac{\alpha_k\eta_k}{(1-\alpha_k)\lambda_k}\sqn{\hx_{k+1} - \hx_k}
  -\tfrac{(1-\alpha_k)\eta_k}{4\alpha_k}\bg_f(\ox_k;\tx_k)
  +\tfrac{\theta\eta_k\eta_{k-1}}{\lambda_k}\bg_f(\ox_{k-1};\tx_{k-1})
  \\&\at{uses \cref{line:alpha,line:eta_H}}\leq
  \tfrac{1}{2}\sqn{x_k - x}
  +H_{k-1}\brr{f(\ox_k)-f(x)}
  -H_k\brr{f(\ox_{k+1}) - f(x)}
  -\tfrac{\eta_k}{4(1+\gamma)}\bg_f(\ox_k;\tx_k)
  \\&
  +\tfrac{(1+\theta)^2}{2\theta}\brr*{\tfrac{\eta_k}{\eta_{k-1}} - \tfrac{\theta + 2\theta^2}{(1+\theta)^2}}\sqn{x_{k+1} - x_k}
  -\tfrac{\theta\eta_k}{2\eta_{k-1}}\sqn{x_k - x_{k-1}}
  \\&
  +\tfrac{\alpha_k\eta_k}{(1-\alpha_k)\lambda_k}\sqn{\hx_{k+1} - \hx_k}
  -\tfrac{(1-\alpha_k)\eta_k}{4\alpha_k}\bg_f(\ox_k;\tx_k)
  +\tfrac{\theta\eta_k\eta_{k-1}}{\lambda_k}\bg_f(\ox_{k-1};\tx_{k-1})
\end{align*}
where \annotate.
Furthermore, we obtain the following:
\begin{align*}
  \mi{1}\tfrac{1}{2}\sqn{x_{k+1} - x}
  \\&\leq
  \tfrac{1}{2}\sqn{x_k - x}
  +H_{k-1}\brr{f(\ox_k)-f(x)}
  -H_k\brr{f(\ox_{k+1}) - f(x)}
  -\tfrac{\eta_k}{4(1+\gamma)}\bg_f(\ox_k;\tx_k)
  \\&
  +\tfrac{(1+\theta)^2}{2\theta}\brr*{\tfrac{\eta_k}{\eta_{k-1}} - \tfrac{\theta + 2\theta^2}{(1+\theta)^2}}\sqn{x_{k+1} - x_k}
  -\tfrac{\theta\eta_k}{2\eta_{k-1}}\sqn{x_k - x_{k-1}}
  \\&
  +\tfrac{\alpha_k\eta_k}{(1-\alpha_k)\lambda_k}\sqn{\hx_{k+1} - \hx_k}
  -\tfrac{(1-\alpha_k)\lambda_{k+1}}{4\theta\alpha_k\eta_{k+1}}\cdot\tfrac{\theta\eta_{k+1}\eta_k}{\lambda_{k+1}}
  \bg_f(\ox_k;\tx_k)
  +\tfrac{\theta\eta_k\eta_{k-1}}{\lambda_k}\bg_f(\ox_{k-1};\tx_{k-1})
  \\&\at{use \cref{line:eta_H,line:alpha}}\leq
  \tfrac{1}{2}\sqn{x_k - x}
  +H_{k-1}\brr{f(\ox_k)-f(x)}
  -H_k\brr{f(\ox_{k+1}) - f(x)}
  -\tfrac{\eta_k}{4(1+\gamma)}\bg_f(\ox_k;\tx_k)
  \\&
  +\tfrac{(1+\theta)^2}{2\theta}\brr*{\tfrac{\eta_k}{\eta_{k-1}} - \tfrac{\theta + 2\theta^2}{(1+\theta)^2}}\sqn{x_{k+1} - x_k}
  -\tfrac{\theta\eta_k}{2\eta_{k-1}}\sqn{x_k - x_{k-1}}
  \\&
  +\tfrac{\alpha_k\eta_k}{(1-\alpha_k)\lambda_k}\sqn{\hx_{k+1} - \hx_k}
  -\tfrac{1}{4\nu\theta(1+\gamma)}\cdot\tfrac{\theta\eta_{k+1}\eta_k}{\lambda_{k+1}}
  \bg_f(\ox_k;\tx_k)
  +\tfrac{\theta\eta_k\eta_{k-1}}{\lambda_k}\bg_f(\ox_{k-1};\tx_{k-1})
  \\&\at{use \cref{eq:nu_gamma}}=
  \tfrac{1}{2}\sqn{x_k - x}
  +H_{k-1}\brr{f(\ox_k)-f(x)}
  -H_k\brr{f(\ox_{k+1}) - f(x)}
  -\tfrac{\eta_k}{4(1+\gamma)}\bg_f(\ox_k;\tx_k)
  \\&
  +\tfrac{(1+\theta)^2}{2\theta}\brr*{\tfrac{\eta_k}{\eta_{k-1}} - \tfrac{\theta + 2\theta^2}{(1+\theta)^2}}\sqn{x_{k+1} - x_k}
  -\tfrac{\theta\eta_k}{2\eta_{k-1}}\sqn{x_k - x_{k-1}}
  \\&
  +\tfrac{\alpha_k\eta_k}{(1-\alpha_k)\lambda_k}\sqn{\hx_{k+1} - \hx_k}
  -\tfrac{(1+\gamma)}{\gamma}\cdot\tfrac{\theta\eta_{k+1}\eta_k}{\lambda_{k+1}}
  \bg_f(\ox_k;\tx_k)
  +\tfrac{\theta\eta_k\eta_{k-1}}{\lambda_k}\bg_f(\ox_{k-1};\tx_{k-1})
  \\&\leq
  \tfrac{1}{2}\sqn{x_k - x}
  +H_{k-1}\brr{f(\ox_k)-f(x)}
  -H_k\brr{f(\ox_{k+1}) - f(x)}
  -\tfrac{\eta_k}{4(1+\gamma)}\bg_f(\ox_k;\tx_k)
  \\&
  +\tfrac{(1+\theta)^2}{2\theta}\brr*{\tfrac{\eta_k}{\eta_{k-1}} - \tfrac{\theta + 2\theta^2}{(1+\theta)^2}}\sqn{x_{k+1} - x_k}
  -\tfrac{\theta\eta_k}{2\eta_{k-1}}\sqn{x_k - x_{k-1}}
  \\&
  +\tfrac{\alpha_k\eta_k}{(1-\alpha_k)\lambda_k}\sqn{\hx_{k+1} - \hx_k}
  -\tfrac{\theta\eta_{k+1}\eta_k}{\lambda_{k+1}}
  \bg_f(\ox_k;\tx_k)
  +\tfrac{\theta\eta_k\eta_{k-1}}{\lambda_k}\bg_f(\ox_{k-1};\tx_{k-1})
  \\&\at{uses \cref{line:hx} and the inequality $\sqn{a+b}\leq 2\sqn{a} + 2\sqn{b}$}\leq
  \tfrac{1}{2}\sqn{x_k - x}
  +H_{k-1}\brr{f(\ox_k)-f(x)}
  -H_k\brr{f(\ox_{k+1}) - f(x)}
  -\tfrac{\eta_k}{4(1+\gamma)}\bg_f(\ox_k;\tx_k)
  \\&
  +\tfrac{\theta\eta_k\eta_{k-1}}{\lambda_k}\bg_f(\ox_{k-1};\tx_{k-1})
  -\tfrac{\theta\eta_{k+1}\eta_k}{\lambda_{k+1}}
  \bg_f(\ox_k;\tx_k)
  \\&
  +\tfrac{(1+\theta)^2}{2\theta}\brr*{\tfrac{\eta_k}{\eta_{k-1}} - \tfrac{\theta + 2\theta^2}{(1+\theta)^2}}\sqn{x_{k+1} - x_k}
  -\tfrac{\theta\eta_k}{2\eta_{k-1}}\sqn{x_k - x_{k-1}}
  \\&
  +\tfrac{2(1+\theta)^2\alpha_k\eta_k}{(1-\alpha_k)\lambda_k}\sqn{x_{k+1} - x_k}
  +\tfrac{2\theta^2\alpha_k\eta_k}{(1-\alpha_k)\lambda_k}\sqn{x_k - x_{k-1}}
  \\&\leq
  \tfrac{1}{2}\sqn{x_k - x}
  +H_{k-1}\brr{f(\ox_k)-f(x)}
  -H_k\brr{f(\ox_{k+1}) - f(x)}
  -\tfrac{\eta_k}{4(1+\gamma)}\bg_f(\ox_k;\tx_k)
  \\&
  +\tfrac{\theta\eta_k\eta_{k-1}}{\lambda_k}\bg_f(\ox_{k-1};\tx_{k-1})
  -\tfrac{\theta\eta_{k+1}\eta_k}{\lambda_{k+1}}
  \bg_f(\ox_k;\tx_k)
  \\&
  +\tfrac{(1+\theta)^2}{2\theta}\brr*{\tfrac{\eta_k}{\eta_{k-1}} + \tfrac{4\theta\alpha_k\eta_k}{(1-\alpha_k)\lambda_k} - \tfrac{\theta + 2\theta^2}{(1+\theta)^2}}\sqn{x_{k+1} - x_k}
  +\tfrac{2\theta^2\alpha_k\eta_k}{(1-\alpha_k)\lambda_k}\sqn{x_k - x_{k-1}}
  \\&\at{uses the inequality $\eta_k \leq (1+\gamma)\eta_{k-1}$ implied by \cref{line:eta_H}}\leq
  \tfrac{1}{2}\sqn{x_k - x}
  +H_{k-1}\brr{f(\ox_k)-f(x)}
  -H_k\brr{f(\ox_{k+1}) - f(x)}
  -\tfrac{\eta_k}{4(1+\gamma)}\bg_f(\ox_k;\tx_k)
  \\&
  +\tfrac{\theta\eta_k\eta_{k-1}}{\lambda_k}\bg_f(\ox_{k-1};\tx_{k-1})
  -\tfrac{\theta\eta_{k+1}\eta_k}{\lambda_{k+1}}
  \bg_f(\ox_k;\tx_k)
  \\&
  +\tfrac{(1+\theta)^2}{2\theta}\brr*{(1 + \gamma) + \tfrac{4\theta\alpha_k\eta_k}{(1-\alpha_k)\lambda_k} - \tfrac{\theta + 2\theta^2}{(1+\theta)^2}}\sqn{x_{k+1} - x_k}
  +\tfrac{2\theta^2\alpha_k\eta_k}{(1-\alpha_k)\lambda_k}\sqn{x_k - x_{k-1}}
  \\&\at{use \cref{line:eta_H,line:alpha}}\leq
  \tfrac{1}{2}\sqn{x_k - x}
  +H_{k-1}\brr{f(\ox_k)-f(x)}
  -H_k\brr{f(\ox_{k+1}) - f(x)}
  -\tfrac{\eta_k}{4(1+\gamma)}\bg_f(\ox_k;\tx_k)
  \\&
  +\tfrac{\theta\eta_k\eta_{k-1}}{\lambda_k}\bg_f(\ox_{k-1};\tx_{k-1})
  -\tfrac{\theta\eta_{k+1}\eta_k}{\lambda_{k+1}}
  \bg_f(\ox_k;\tx_k)
  \\&
  +\tfrac{(1+\theta)^2}{2\theta}\brr*{(1 + \gamma) + 4\nu\theta(1+\gamma)^2 - \tfrac{\theta + 2\theta^2}{(1+\theta)^2}}\sqn{x_{k+1} - x_k}
  +2\nu\theta^2(1+\gamma)^2\sqn{x_k - x_{k-1}}
  \\&\at{use \cref{eq:nu_gamma}}=
  \tfrac{1}{2}\sqn{x_k - x}
  +H_{k-1}\brr{f(\ox_k)-f(x)}
  -H_k\brr{f(\ox_{k+1}) - f(x)}
  -\tfrac{\eta_k}{4(1+\gamma)}\bg_f(\ox_k;\tx_k)
  \\&
  +\tfrac{\theta\eta_k\eta_{k-1}}{\lambda_k}\bg_f(\ox_{k-1};\tx_{k-1})
  -\tfrac{\theta\eta_{k+1}\eta_k}{\lambda_{k+1}}
  \bg_f(\ox_k;\tx_k)
  \\&
  +\tfrac{(1+\theta)^2}{2\theta}\brr*{1 + 2\gamma - \tfrac{\theta}{(1+\theta)} - \tfrac{\theta^2}{(1+\theta)^2}}\sqn{x_{k+1} - x_k}
  +\tfrac{\gamma\theta}{2}\sqn{x_k - x_{k-1}}
  \\&\at{use \cref{eq:nu_gamma}}\leq
  \tfrac{1}{2}\sqn{x_k - x}
  +H_{k-1}\brr{f(\ox_k)-f(x)}
  -H_k\brr{f(\ox_{k+1}) - f(x)}
  -\tfrac{\eta_k}{4(1+\gamma)}\bg_f(\ox_k;\tx_k)
  \\&
  +\tfrac{\theta\eta_k\eta_{k-1}}{\lambda_k}\bg_f(\ox_{k-1};\tx_{k-1})
  -\tfrac{\theta\eta_{k+1}\eta_k}{\lambda_{k+1}}
  \bg_f(\ox_k;\tx_k)
  -\gamma\theta\sqn{x_{k+1} - x_k}
  +\tfrac{\gamma\theta}{2}\sqn{x_k - x_{k-1}}
  \\&\at{desc}=
  \tfrac{1}{2}\sqn{x_k - x}
  +H_{k-1}\brr{f(\ox_k)-f(x)}
  -H_k\brr{f(\ox_{k+1}) - f(x)}
  \\&
  +\tfrac{\theta\eta_k\eta_{k-1}}{\lambda_k}\bg_f(\ox_{k-1};\tx_{k-1})
  -\tfrac{\theta\eta_{k+1}\eta_k}{\lambda_{k+1}}
  \bg_f(\ox_k;\tx_k)
  -\tfrac{\gamma\theta}{2}\sqn{x_{k+1} - x_k}
  +\tfrac{\gamma\theta}{2}\sqn{x_k - x_{k-1}}
  \\&
  -\tfrac{\gamma\theta\eta_k^2}{2}\sqn{\nabla f(\tx_k)}
  -\tfrac{\eta_k}{4(1+\gamma)}\bg_f(\ox_k;\tx_k)m
\end{align*}
where \annotate. It remains to use the definition $\Psi_k(x)$ in \cref{eq:Psi}.\qed

\proofsubsection{cor:main}

We can upper-bound $\tfrac{1}{2}\sqn{x_K - x} + H_{K-1}\brr{f(\ox_K)-f(x)}$ as follows:
\begin{align*}
  \mi{3}
  \tfrac{1}{2}\sqn{x_K - x} + H_{K-1}\brr{f(\ox_K)-f(x)}
  \\&\at{use \cref{eq:Psi}}\leq
  \Psi_K(x)
  \at{uses \Cref{thm:main}}\leq
  \Psi_1(x)
  \\&\at{use \cref{eq:Psi}}\leq
  \tfrac{1}{2}\sqn{x_1 - x}
  +H_{0}\brr{f(\ox_1)-f(x)}
  +\tfrac{\theta\eta_k\eta_{k-1}}{\lambda_k}\bg_f(\ox_{0};\tx_{0})
  +\tfrac{\gamma\theta}{2}\sqn{x_1 - x_0}
  \\&\at{uses \cref{line:init,line:ox}}=
  \tfrac{1}{2}\sqn{x_1 - x}
  +\eta_0\brr{f(x_0)-f(x)}
  +\tfrac{\gamma\theta}{2}\sqn{x_1 - x_0}
  \\&\at{use \cref{line:x,line:init}}=
  \tfrac{1}{2}\sqn{x_0 - x}
  +\tfrac{(1+\gamma\theta)\eta_0^2}{2}\sqn{\nabla f(x_0)}
  -\eta_0\<\nabla f(x_0),x_0 - x>
  +\eta_0\brr{f(x_0)-f(x)}
  \\&=
  \tfrac{1}{2}\sqn{x_0 - x}
  +\tfrac{(1+\gamma\theta)\eta_0^2}{2}\sqn{\nabla f(x_0)}
  -\eta_0\bg_f(x,x_0)
  \\&\at{uses the convexity of $f(x)$}\leq
  \tfrac{1}{2}\sqn{x_0 - x}
  +\tfrac{(1+\gamma\theta)\eta_0^2}{2}\sqn{\nabla f(x_0)},
\end{align*}
where \annotate.\qed

\newpage

\section[Proofs for \crtCref{sec:L}]{Proofs for \Cref{sec:L}}

\proofsubsection{lem:H}

For $k=0$ we have $H_k = H_{k-1}$. The inequality $H_k \geq H_{k-1}$ is obvious for $k \geq 1$ due to \cref{line:eta_H}. Furthermore, we can upper-bound $H_k$ for $k \geq 1$ as follows:
\begin{align*}
  H_k
  &\at{use \cref{line:eta_H}}=
  H_{k-1} +\eta_k
  \at{use \cref{line:eta_H}}\leq
  H_{k-1} + (1+\gamma)\eta_{k-1}
  \at{use \cref{line:eta_H}}\leq
  (2+\gamma)H_{k-1},
\end{align*}
where \annotate.\qed

\proofsubsection{thm:L}

We prove the statement of \Cref{thm:L} by induction over $K$. The base case, $K=m$, is obvious. Now, we assume that \cref{eq:H_L} holds for all $k \in \brf{0,\dots,K}$, where $K\geq m$, and prove \cref{eq:H_L} for $k = K+1$. We consider the following cases:
\begin{enumerate}[label=\textbf{Case \arabic*.}]
  \item  $(1+\gamma) \eta_K\eta_{K-1} \geq \nu H_{K-1}\lambda_{K+1}$.
  \item $(1+\gamma) \eta_K\eta_{K-1} < \nu H_{K-1}\lambda_{K+1}$.
    \begin{enumerate}[label=\textbf{Case 2\alph*.}]
      \item $\eta_K = \eta_0 (1+\gamma)^K$.
      \item $\eta_K < \eta_0 (1+\gamma)^K$.
    \end{enumerate}
\end{enumerate}

\textbf{Case 1.} Here we have $\eta_{K+1} = \frac{\nu H_{K-1}\lambda_{k+1}}{\eta_{K-1}}$ and obtain the following inequality:
\begin{align*}
  \sqrt{H_{K+1}}
  &=
  \tfrac{1}{\sqrt{H_{K+1}} + \sqrt{H_{K-2}}} \brr{H_{K+1} - H_{K-2}} + \sqrt{H_{K-2}}
  \\&\at{uses the inequality $H_{K-2} \leq H_{K+1}$ for $K \geq m \geq 2$}\geq
  \tfrac{1}{2\sqrt{H_{K+1}}} \brr{H_{K+1} - H_{K-2}} + \sqrt{H_{K-2}}
  \\&\at{uses \cref{line:eta_H}}\geq
  \tfrac{1}{2\sqrt{H_{K+1}}} \brr{\eta_{K+1} + \eta_{K-1}} + \sqrt{H_{K-2}}
  \\&\at{uses Young's inequality}\geq
  \sqrt{\tfrac{\eta_{K+1}\eta_{K-1}}{H_{K+1}}} + \sqrt{H_{K-2}}
  \\&\at{uses the assumption $\eta_{K+1} = \frac{\nu H_{K-1}\lambda_{k+1}}{\eta_{K-1}}$ (Case 1) and \Cref{lem:Lambda_L}}\geq
  \sqrt{\tfrac{\nu H_{K-1}}{LH_{K+1}}} + \sqrt{H_{K-2}}
  \\&\at{uses \Cref{lem:H}}\geq
  \tfrac{\sqrt{\nu}}{(2+\gamma)\sqrt{L}} + \sqrt{H_{K-2}}
  \\&\at{uses the induction hypothesis in \cref{eq:H_L}}\geq
  \tfrac{\sqrt{\nu}}{(2+\gamma)\sqrt{L}} + \tfrac{c}{\sqrt{L}}\cdot(K - 2 - m)
  \\&\at{uses the inequality $\frac{3c}{\sqrt{L}} \leq \frac{\sqrt{\nu}}{(2+\gamma)\sqrt{L}}$ implied by \cref{eq:cm_L}}\geq
  \tfrac{c}{\sqrt{L}}\cdot(K + 1 - m),
\end{align*}
where \annotate. This proves Case~1.

\textbf{Case 2a.} In this case, it is easy to verify that $\eta_k = \eta_0(1+\gamma)^k$ for all $k \in \brf{0,\dots,K+1}$. Hence, we can lower-bound $H_{K+1}$ as follows:
\begin{align}
  H_{K+1}
  &\at{uses \cref{line:eta_H} and the relation $\eta_k = \eta_0(1+\gamma)^k$ above}=\notag
  \tsum_{k=0}^{K+1}\eta_0(1+\gamma)^k
  =
  \tfrac{\eta_0}{\gamma}\brr*{(1+\gamma)^{K+2} - 1}
  \\&=\notag
  \tfrac{\eta_0(1+\gamma)^m}{\gamma}\brr*{(1+\gamma)^{K+2-m} - \tfrac{1}{(1+\gamma)^m}}
  \geq
  \tfrac{\eta_0(1+\gamma)^m}{\gamma}\brr*{\brr*{(1+\gamma)^{(K-m+2)/2}}^2 - 1}
  \\&\at{uses the inequality $K\geq m$ and Bernoulli's inequality}\geq
  \label{aux:L:2a}
  \tfrac{\eta_0(1+\gamma)^m}{\gamma}\brr*{\brr*{1 + \tfrac{\gamma}{2}(K-m+2)}^2 - 1}
  \geq
  \tfrac{\eta_0\gamma(1+\gamma)^m}{4}\brr*{K-m+2}^2
  \\&\at{uses the inequality $(1+\gamma)^m \geq \frac{4c^2}{\gamma\eta_0 L}$ implied by \cref{eq:cm_L}}\geq\notag
  \tfrac{c^2}{L}(K+2-m)^2
  \geq
  \tfrac{c^2}{L}(K+1-m)^2,
\end{align}
where \annotate. This proves Case~2a.

\textbf{Case 2b.} In this case, there exists $l \in \brf{1,...,K}$ such that $\eta_{l} < (1+\gamma)\eta_{l-1}$. We choose the largest such index $l \in \brf{1,...,K}$. This implies the following relations:
\begin{equation}
  \eta_l = \tfrac{\nu H_{l-2}\lambda_l}{\eta_{l-2}}
  \quad\text{and}\quad
  \eta_{k} = (1+\gamma)^{k-l}\eta_l
  \;\;\text{for all}\;\;
  k \in \brf{l,\dots,K+1}.
\end{equation}
Hence, we can upper-bound $H_{K+1}$ as follows:
\begin{equation}\label{aux:L:H}
  \begin{aligned}
    H_{K+1}
    &\at{uses the above relation}=
    H_{l-1} + \tsum_{i=0}^{K-l+1}\eta_l(1+\gamma)^{i}
    =
    H_{l-1} + \tfrac{\eta_l}{\gamma}\brr{(1+\gamma)^{K-l+2} - 1}
    \\&\at{uses the inequalities $\eta_l \leq (1+\gamma)\eta_{l-1}$ and $\eta_{l-1} \leq H_{l-1}$ implied by \cref{line:eta_H}}\leq
    \brr*{1 + \tfrac{(1+\gamma)^{K-l+3} - (1+\gamma)}{\gamma}} H_{l-1}
    \leq
    \tfrac{(1+\gamma)^{K-l+3}}{\gamma} H_{l-1}
  \end{aligned}
\end{equation}
where \annotate.
Furthermore, we can lower-bound $\eta_l$ as follows:
\begin{equation}\label{aux:L:eta}
  \eta_l
  \at{uses the relation above}=
  \tfrac{\nu H_{l-2}\lambda_l}{\eta_{l-2}}
  \at{uses \cref{line:eta_H,line:init}, where we define $H_{-2} = H_0$}\geq
  \tfrac{\nu H_{l-2}\lambda_l}{H_{l} - H_{l-3}}
  \at{uses \Cref{lem:Lambda_L}}\geq
  \tfrac{\nu H_{l-2}}{L(H_{l} - H_{l-3})},
\end{equation}
where \annotate.
Next, we obtain the following inequality:
\begin{align}
  \brr*{\sqrt{H_{K+1}} - \sqrt{H_{l-3}}}^2
  &\at{uses the inequalities $H_{l-3} \leq H_{l-1}$ and $H_l\leq H_{K+1}$}\geq\notag
  \brr*{\sqrt{H_{K+1}} - \sqrt{H_{l-1}}}\brr*{\sqrt{H_l} - \sqrt{H_{l-3}}}
  \\&=\notag
  \tfrac{H_{K+1} - H_{l-1}}{\sqrt{H_{K+1}} + \sqrt{H_{l-1}}}\cdot\tfrac{H_l - H_{l-3}}{\sqrt{H_l} + \sqrt{H_{l-3}}}
  \\&\at{uses the inequalities $H_{K+1} \geq H_{l-1}$ and $H_l \geq H_{l-3}$}\geq\notag
  \tfrac{H_{K+1} - H_{l-1}}{2\sqrt{H_{K+1}}}\cdot\tfrac{H_l - H_{l-3}}{2\sqrt{H_l}}
  \\&\at{use \cref{aux:L:H}}=\notag
  \tfrac{\eta_l\brr{(1+\gamma)^{K-l+2} - 1}}{2\gamma\sqrt{H_{K+1}}}\cdot\tfrac{H_l - H_{l-3}}{2\sqrt{H_l}}
  \\&\at{use \cref{aux:L:eta}}\geq\notag
  \tfrac{\nu H_{l-2}\brr{(1+\gamma)^{K-l+2} - 1}}{4\gamma L\sqrt{H_{K+1}H_l}}
  \\&\at{use \cref{aux:L:H}}\geq\notag
  \tfrac{\nu H_{l-2}}{4 L\sqrt{\gamma(1+\gamma) H_{l-1}H_l}}\cdot \tfrac{(1+\gamma)^{K-l+2} - 1}{\sqrt{(1+\gamma)^{K-l+2}}}
  \\&\geq\notag
  \tfrac{\nu H_{l-2}}{4 L\sqrt{\gamma(1+\gamma) H_{l-1}H_l}}\cdot\brr{(1+\gamma)^{(K-l)/2+1} - 1}
  \\&\at{uses \Cref{lem:H}}\geq\label{aux:L:2b}
  \tfrac{\nu}{4 L\sqrt{\gamma(1+\gamma)(2+\gamma)^3}}\cdot\brr{(1+\gamma)^{(K-l)/2+1} - 1},
\end{align}
where \annotate. Next, we take the square root of both sides and obtain the following:
\begin{align}
  \sqrt{H_{K+1}} - \sqrt{H_{l-3}}
  &\geq\notag
  \tfrac{\sqrt{\nu}}{2\sqrt[4]{\gamma(1+\gamma)(2+\gamma)^3L^2}}\cdot\sqrt{(1+\gamma)^{(K-l)/2+1} - 1}
  \\&\at{uses the inequality $\sqrt{a} \geq \sqrt{a+b} - \sqrt{b}$ for $a,b\geq 0$}\geq\notag
  \tfrac{\sqrt{\nu}}{2\sqrt[4]{\gamma(1+\gamma)(2+\gamma)^3L^2}}\cdot\brr*{(1+\gamma)^{(K-l+2)/4} - 1}
  \\&=\notag
  \tfrac{\sqrt{\nu}}{2\sqrt[4]{\gamma(1+\gamma)^5(2+\gamma)^3L^2}}\cdot\brr*{(1+\gamma)^{(K-l)/4+3/2} - \brr{1+\gamma}}
  \\&\at{uses Bernoulli's inequality}\geq\notag
  \tfrac{\sqrt{\nu}\gamma}{2\sqrt[4]{\gamma(1+\gamma)^5(2+\gamma)^3L^2}}\cdot\brr*{\tfrac{1}{4}(K-l)+\tfrac{1}{2}}
  \\&\at{uses the inequality $K \geq l$}\geq\notag
  \tfrac{\sqrt{\nu}\gamma}{2\sqrt[4]{\gamma(1+\gamma)^5(2+\gamma)^3L^2}}\cdot\brr*{\tfrac{1}{8}(K-l)+\tfrac{1}{2}}
  \\&=\label{aux:L:2b_1}
  \tfrac{\sqrt{\nu}\gamma}{16\sqrt[4]{\gamma(1+\gamma)^5(2+\gamma)^3L^2}}\cdot\brr*{K-l+4}
  \\&\at{uses \cref{eq:cm_L}}\geq\notag
  \tfrac{c}{\sqrt{L}}\brr*{K-l+4},
\end{align}
where \annotate. Finally, we obtain the following:
\begin{align*}
  \sqrt{H_{K+1}} \geq \sqrt{H_{l-3}} + \tfrac{c}{\sqrt{L}}\brr*{K-l+4}
  \at{uses the induction hypothesis in \cref{eq:H_L} for $k=l-3$, which also holds for $l=1$ due to the definition $H_{-2} = H_0$}\geq
  \tfrac{c}{\sqrt{L}}\brr*{K+1 - m},
\end{align*}
where \annotate. This proves Case~2b.\qed

\newpage

\section[Proofs for \crtCref{sec:L01}]{Proofs for \Cref{sec:L01}}

\proofsubsection{lem:D}

\textbf{Step~1.}
We can upper-bound $\sum_{i=1}^K\brr*{\eta_i^2\sqn{\nabla f(\tx_i)} + \eta_i\bg_f(\ox_i;\tx_i)}$ as follows:
\begin{align*}
  \mi{2}
  \tsum_{i=1}^K\brr*{\eta_i^2\sqn{\nabla f(\tx_i)} + \eta_i\bg_f(\ox_i;\tx_i)}
  \\&\leq
  \max\brf*{\tfrac{2}{\gamma\theta},4(1+\gamma)}\tsum_{i=1}^K\brr*{\tfrac{\gamma\theta}{2}\eta_i^2\sqn{\nabla f(\tx_i)} + \tfrac{1}{4(1+\gamma)}\eta_i\bg_f(\ox_i;\tx_i)}
  \\&\at{uses \Cref{thm:main}}\leq
  \max\brf*{\tfrac{2}{\gamma\theta},4(1+\gamma)}\tsum_{i=1}^K\brr*{
    \Psi_i(x^*) - \Psi_{i+1}(x^*)
  }
  \\&\at{uses the fact that $\Psi_1(x^*)\geq 0$, which is implied by \cref{eq:Psi}}\leq
  \max\brf*{\tfrac{2}{\gamma\theta},4(1+\gamma)}\Psi_1(x^*)
  \\&\at{is obtained similarly to the proof of \Cref{cor:main} in \Cref{proof:cor:main}}\leq
  \max\brf*{\tfrac{1}{\gamma\theta},2(1+\gamma)}
  \brr*{
    \sqn{x_0 - x^*}
  +(1+\gamma\theta)\eta_0^2\sqn{\nabla f(x_0)}},
\end{align*}
where \annotate.

\textbf{Step~2.}
Next, we prove the following inequality for all $k \in \brf{0,\dots,K}$ by induction:
\begin{equation}\label{aux:dist}
  \max\brf*{\norm{\tx_k - x^*},\norm{\ox_k - x^*}} \leq
  (1+2\theta)\max_{k=0,\dots,K}\norm{x_k - x^*}.
\end{equation}
The base case $k=0$ is obvious due to \cref{line:init}. For $k\geq 1$ we can upper-bound $\norm{\ox_k - x^*}$ as follows:
\begin{align*}
  \norm{\ox_k - x^*}
  &\at{uses \cref{line:ox}, the triangle inequality, and \Cref{lem:alpha_beta}}=
  \beta_{k-1}\norm{\tx_{k-1} - x^*} + (1-\beta_{k-1})\norm{\ox_{k-1} - x^*}
  \\&\at{uses the induction hypothesis in \cref{aux:dist} for $k-1$}\leq
  (1+2\theta)\max_{k=0,\dots,K}\norm{x_k - x^*},
\end{align*}
where \annotate. Next, we can upper-bound $\norm{\tx_k - x^*}$ as follows:
\begin{align*}
  \norm{\tx_k - x^*}
  &\at{uses \cref{line:tx} and the triangle inequality}\leq
  \alpha_k\norm{\hx_k - x^*} + (1-\alpha_k)\norm{\ox_k - x^*}
  \\&\at{uses the inequality obtained above}\leq
  \alpha_k\norm{\hx_k - x^*} + (1-\alpha_k)(1+2\theta)\max_{k=0,\dots,K}\norm{x_k - x^*}
  \\&\at{uses \cref{line:hx} and the triangle inequality}\leq
  \alpha_k(1+\theta)\norm{x_k - x^*} + \alpha_k\theta\norm{x_{k-1} - x^*} + (1-\alpha_k)(1+2\theta)\max_{k=0,\dots,K}\norm{x_k - x^*}
  \\&\leq
  (1+2\theta)\max_{k=0,\dots,K}\norm{x_k - x^*},
\end{align*}
where \annotate. This proves \cref{aux:dist}.
Next, using \cref{aux:dist} and \Cref{cor:main}, we obtain the following inequality for $k \in \brf{0,\dots,K}$:
\begin{align*}
  \max\brf*{\sqn{\tx_k - x^*},\sqn{\ox_k - x^*}}
  \leq
  (1+2\theta)^2
  \brr*{
    \sqn{x_0 - x^*}
  +(1+\gamma\theta)\eta_0^2\sqn{\nabla f(x_0)}}.
\end{align*}
Combining this with the inequality obtained in Step~1 concludes the proof.\qed

\proofsubsection{lem:lambda_min}

We can lower-bound $\lambda_k$ as follows:
\begin{align*}
  \lambda_k
  &\at{uses \cref{line:lambda}}=
  \min\brf*{\Lambda(\ox_k,\tx_{k-1}),\Lambda(\ox_k,\tx_k)}
  \\&\at{uses \cref{eq:Lambda}}=
  \min\brf*{\tfrac{2\bg_f(\ox_k,\tx_{k-1})}{\sqn{\nabla f(\ox_k) - \nabla f(\tx_{k-1})}},\tfrac{2\bg_f(\ox_k,\tx_k)}{\sqn{\nabla f(\ox_k) - \nabla f(\tx_k)}}}
  \\&\at{uses Corollary~2.8 of \citet{vankov2024optimizing}}\geq
  \frac{2}{2(L_0 + L_1\norm{\nabla f(\ox_k)}) + L_1\max\brf*{\norm{\nabla f(\ox_k)-\nabla f(\tx_k)},\norm{\nabla f(\ox_k)-\nabla f(\tx_{k-1})}}}
  \\&\at{use \cref{eq:L01}}\geq
  \frac{2}{\brr{L_0 + L_1\norm{\nabla f(\ox_k)}}\brr{\exp\brr*{L_1\max\brf*{\norm{\ox_k-\tx_k},\norm{\ox_k-\tx_{k-1}}}} + 1}}
  \\&\geq
  \frac{2}{\brr{L_0 + L_1\norm{\nabla f(\ox_k)}}\brr{\exp\brr*{L_1\brr{\norm{\ox_k-x^*} + \max\brf*{\norm{\tx_k - x^*},\norm{\tx_{k-1}-x^*}}}} + 1}}
  \\&\at{use \Cref{lem:D}}\geq
  \frac{2}{\brr{L_0 + L_1\norm{\nabla f(\ox_k)}}\brr{1+\exp\brr*{2L_1\cD}}}
  \\&\at{use \cref{eq:L01}}\geq
  \frac{2}{L_0\exp\brr{L_1\norm{\ox_k-x^*}}\brr{1+\exp\brr*{2L_1\cD}}}
  \\&\at{use \Cref{lem:D}}\geq
  \frac{2}{L_0\exp\brr{L_1\cD}\brr{1+\exp\brr*{2L_1\cD}}}
  \geq
  \frac{1}{L_0\exp\brr{3L_1\cD}},
\end{align*}
where \annotate.\qed

\proofsubsection{lem:lambda_k}

We can lower-bound $\bg_f(x,z)$ for all $x,z\in\sX$ as follows:
\begin{align*}
  \bg_f(x,z)
  &\at{uses Corollary~2.8 of \citet{vankov2024optimizing}}\geq
  \frac{\sqn{\nabla f(x) - \nabla f(z)}}{2(L_0 + L_1\norm{\nabla f(x)}) + L_1\norm{\nabla f(x)-\nabla f(z)}}
  \\&\at{uses the triangle inequality}\geq
  \frac{\sqn{\nabla f(x) - \nabla f(z)}}{2(L_0 + L_1\norm{\nabla f(z)}) + 3L_1\norm{\nabla f(x)-\nabla f(z)}},
\end{align*}
where \annotate. Hence, we obtain the following inequality, which is quadratic in $\norm{\nabla f(x)-\nabla f(z)}$:
\begin{align*}
  \sqn{\nabla f(x) - \nabla f(z)} - 3L_1\bg_f(x,z) \cdot \norm{\nabla f(x)-\nabla f(z)}
  - 2(L_0 + L_1\norm{\nabla f(z)})\bg_f(x,z) \leq 0.
\end{align*}
Solving this inequality with respect to $\norm{\nabla f(x)-\nabla f(z)}$ gives the following
\begin{align*}
  \norm{\nabla f(x)-\nabla f(z)}
  &\leq
  \frac{3L_1\bg_f(x,z) + \sqrt{\brr{3L_1\bg_f(x,z)}^2 + 8(L_0 + L_1\norm{\nabla f(z)})\bg_f(x,z)}}{2}
  \\&\at{uses the inequality $\sqrt{a+b}\leq \sqrt{a}+\sqrt{b}$}\leq
  3L_1\bg_f(x,z) + \sqrt{2(L_0 + L_1\norm{\nabla f(z)})\bg_f(x,z)}
  \\&\at{uses Young's inequality}\leq
  \tfrac{9}{2}L_1\bg_f(x,z) + \tfrac{1}{3}(L_0/L_1 + \norm{\nabla f(z)}),
\end{align*}
where \annotate. Combining this with the inequalities above gives the following:
\begin{align*}
  \Lambda(x,z)
  &\at{uses \cref{eq:Lambda}}=
  \tfrac{2\bg_f(x,z)}{\sqn{\nabla f(x) - \nabla f(z)}}
  \\&\at{use the inequalities obtained above}\geq
  \frac{2}{2(L_0 + L_1\norm{\nabla f(z)}) + 3L_1\norm{\nabla f(x)-\nabla f(z)}}
  \\&\at{use the inequalities obtained above}\geq
  \frac{4}{6(L_0 + L_1\norm{\nabla f(z)}) + 27L_1^2\bg_f(x,z)}
  \\&\geq
  \frac{4}{9\max\brf*{2L_0, 2L_1\norm{\nabla f(z)},9L_1^2\bg_f(x,z)}}
\end{align*}
where \annotate. Using this, we can lower-bound $\lambda_k$ as follows:
\begin{align*}
  \lambda_k
  &\at{uses \cref{line:lambda}}\geq
  \min\brf*{\Lambda(\ox_k,\tx_{k-1}),\Lambda(\ox_k,\tx_k)}
  \\&\at{uses the inequality obtained above}\geq
  \frac{4}{9\max\brf*{2L_0, 2L_1\norm{\nabla f(\tx_k)},2L_1\norm{\nabla f(\tx_{k-1})},9L_1^2\bg_f(\ox_k,\tx_{k}),9L_1^2\bg_f(\ox_k,\tx_{k-1})}}
  \\&\at{uses \Cref{lem:beta_f}}\geq
  \frac{4}{9\max\brf*{2L_0, 2L_1\norm{\nabla f(\tx_k)},2L_1\norm{\nabla f(\tx_{k-1})},9L_1^2\bg_f(\ox_k,\tx_{k}),9L_1^2\bg_f(\ox_{k-1},\tx_{k-1})}},
\end{align*}
where \annotate.\qed

\proofsubsection{lem:Tau}

Using the definition of $\cT_4(k)$ in \cref{eq:Tau}, one can verify that the following inclusion holds:
\begin{equation}
  \cT_4(k)\subset \brf*{1 \leq i \leq k : i - \max \brf*{j : j \in \cT_2(i)\cup\brf{0}} \leq m},
\end{equation}
which implies the first desired inequality. Next, using \Cref{lem:lambda_k} and the definition of $\cT_2(k)$ in \cref{eq:Tau}, one can veryfy that the following inequality holds for all $i \in \cT_2(k)$:
\begin{equation}\label{aux:L01:lambda}
  \lambda_i \geq \frac{4}{9\max\brf*{2L_1\norm{\nabla f(\tx_i)},2L_1\norm{\nabla f(\tx_{i-1})},9L_1^2\bg_f(\ox_i,\tx_i),9L_1^2\bg_f(\ox_{i-1},\tx_{i-1})}}.
\end{equation}
Hence, we obtain the following:
\begin{align*}
  \mi{0.5}2\cD^2
  \\&\at{uses \Cref{lem:D}}\geq
  2\tsum_{i=1}^k\brr*{\eta_i^2\sqn{\nabla f(\tx_i)} + \eta_i\bg_f(\ox_i;\tx_i)}
  \\&=
  \tsum_{i=1}^k\brr*{\eta_i^2\sqn{\nabla f(\tx_i)} + \eta_i\bg_f(\ox_i;\tx_i)}
  \\&
  +\tsum_{i=2}^{k+1}\brr*{\eta_{i-1}^2\sqn{\nabla f(\tx_{i-1})} + \eta_{i-1}\bg_f(\ox_{i-1};\tx_{i-1})}
  \\&\at{uses \cref{line:eta_H}}\geq
  \tsum_{i=1}^k\brr*{\eta_i^2\sqn{\nabla f(\tx_i)} + \eta_i\bg_f(\ox_i;\tx_i)}
  \\&
  +\tsum_{i=2}^{k+1}\brr*{\tfrac{1}{(1+\gamma)^2}\eta_{i}^2\sqn{\nabla f(\tx_{i-1})} + \tfrac{1}{(1+\gamma)}\eta_{i}\bg_f(\ox_{i-1};\tx_{i-1})}
  \\&\at{uses the fact that $\gamma > 0$}\geq
  \tsum_{i=2}^k
  \brr*{
    \tfrac{\eta_i^2}{(1+\gamma)^2}\brr*{\sqn{\nabla f(\tx_i)} + \sqn{\nabla f(\tx_{i-1})}}
    +
    \tfrac{\eta_i}{(1+\gamma)^2}\brr*{\bg_f(\ox_i;\tx_i) + \bg_f(\ox_{i-1};\tx_{i-1})}
  }
  \\&\at{uses the definition $\cT_2'(k) = \cT_2(k)\setminus\brf{1}$}\geq
  \tsum_{i\in \cT_2'(k)}
  \brr*{
    \tfrac{\eta_i^2}{(1+\gamma)^2}\brr*{\sqn{\nabla f(\tx_i)} + \sqn{\nabla f(\tx_{i-1})}}
    +
    \tfrac{\eta_i}{(1+\gamma)^2}\brr*{\bg_f(\ox_i;\tx_i) + \bg_f(\ox_{i-1};\tx_{i-1})}
  }
  \\&\at{uses the definition of $\cT_2(k)$ in \cref{eq:Tau} and the fact that $H_{i-2}\geq \eta_{i-2}$}\geq
  \tsum_{i\in \cT_2'(k)}
  \brr*{
    \tfrac{\nu^2\lambda_i^2}{(1+\gamma)^2}\brr*{\sqn{\nabla f(\tx_i)} + \sqn{\nabla f(\tx_{i-1})}}
    +
    \tfrac{\nu\lambda_i}{(1+\gamma)^2}\brr*{\bg_f(\ox_i;\tx_i) + \bg_f(\ox_{i-1};\tx_{i-1})}
  },
\end{align*}
where \annotate. Furthermore, using \cref{aux:L01:lambda}, we can show that the following inequality holds for all $i \in \cT_2'(k)$:
\begin{equation}
  \max\brf*{
    \lambda_i^2\sqn{\nabla f(\tx_i)},
    \lambda_i^2\sqn{\nabla f(\tx_{i-1})},
    \lambda_i\bg_f(\ox_i;\tx_i),
    \lambda_i\bg_f(\ox_{i-1};\tx_{i-1}),
  } \geq \tfrac{4}{81L_1^2}.
\end{equation}
Hence, we obtain the following:
\begin{align*}
  2\cD^2
  &\geq
  \tsum_{i\in \cT_2'(k)}
  \brr*{
    \tfrac{\min\brf{\nu,\nu^2}}{(1+\gamma)^2}\cdot
    \tfrac{4}{81L_1^2}
  }
  \at{uses the fact that $\gamma > 0$}\geq
  \tfrac{4\min\brf{\nu,\nu^2}}{81(1+\gamma)^2L_1^2}\cdot \abs{\cT_2'(k)},
\end{align*}
where \annotate. It remains to use the inequality $\abs{\cT_2(k)} \leq 1 + \abs{\cT_2'(k)}$.\qed

\proofsubsection{thm:L01}

We prove the statement of \Cref{thm:L01} by induction over $K$. The base case, $K=1$, is obvious. Now, we assume that \cref{eq:H_L01} holds for all $k \in \brf{0,\dots,K}$, where $K\geq 1$, and prove the inequality in \cref{eq:H_L01} for $k = K+1$. We consider the following cases:
\begin{enumerate}[label=\textbf{Case \arabic*.}]
  \item $K+1\in\cT_2(K+1)$ or $K+1\in\cT_4(K+1)$.
  \item $K+1 \in \cT_1(K+1)$.
  \item $K+1 \in \cT_3(K+1)$.
    \begin{enumerate}[label=\textbf{Case 3\alph*.}]
      \item $l(K+1) = 0$.
      \item $l(K+1) \in \cT_1(K+1)$.
      \item $l(K+1) \in \cT_2(K+1)$.
    \end{enumerate}
\end{enumerate}

\textbf{Case 1.} In this case, we can lower-bound $\sqrt{H_{K+1}}$ as follows:
\begin{align*}
  \sqrt{H_{K+1}}
  &\at{uses \Cref{lem:H}}\geq
  \sqrt{H_K}
  \\&\at{uses the induction hypothesis in \cref{eq:H_L01}}\geq
  \tfrac{c}{\sqrt{L_0}} \cdot \brr*{K - \abs{\cT_2(K)} - \abs{\cT_4(K)} - 1}
  \\&\at{uses the assumpttion $K+1 \in \cT_2(K+1)\cup\cT_4(K+1)$ and the definition of $\cT_2(K+1)$ and $\cT_4(K+1)$ in \cref{eq:Tau}}\geq
  \tfrac{c}{\sqrt{L_0}} \cdot \brr*{K+1 - \abs{\cT_2(K+1)} - \abs{\cT_4(K+1)} - 1},
\end{align*}
where \annotate. This proves Case~1.

\textbf{Case 2.}
In this case, we can lower-bound $\sqrt{H_{K+1}}$ as follows:
\begin{align*}
  \sqrt{H_{K+1}}
  &=
  \tfrac{1}{\sqrt{H_{K+1}} + \sqrt{H_{K-2}}}\brr{H_{K+1} - H_{K-2}}  + \sqrt{H_{K-2}}
  \\&\at{use \Cref{lem:H}}\geq
  \tfrac{1}{2\sqrt{H_{K+1}}}\brr{H_{K+1} - H_{K-2}}  + \sqrt{H_{K-2}}
  \\&\at{uses \cref{line:eta_H}}\geq
  \tfrac{1}{2\sqrt{H_{K+1}}}\brr{\eta_{K+1} + \eta_{K-1}}  + \sqrt{H_{K-2}}
  \\&\at{uses Young's inequality}\geq
  \sqrt{\tfrac{\eta_{K+1}\eta_{K-1}}{H_{K+1}}}  + \sqrt{H_{K-2}}
  \\&\at{use the assumption $K+1 \in \cT_1(K+1)$ and the definition of $\cT_1(K+1)$ in \cref{eq:Tau}}\geq
  \sqrt{\tfrac{\nu \lambda_{K+1}H_{K-1}}{H_{K+1}}}  + \sqrt{H_{K-2}}
  \\&\at{use the assumption $K+1 \in \cT_1(K+1)$ and the definition of $\cT_1(K+1)$ in \cref{eq:Tau}}\geq
  \sqrt{\tfrac{2\nu H_{K-1}}{9 L_0 H_{K+1}}}  + \sqrt{H_{K-2}}
  \\&\at{use \Cref{lem:H}}\geq
  \tfrac{\sqrt{2\nu}}{3(2+\gamma)\sqrt{L_0}} + \sqrt{H_{K-2}}
  \\&\at{uses the definition of $c$ in \cref{eq:cm_L01}}\geq
  \tfrac{3c}{\sqrt{L_0}} + \sqrt{H_{K-2}}
  \\&\at{uses the induction hypothesis in \cref{eq:H_L01}}\geq
  \tfrac{3c}{\sqrt{L_0}} + \tfrac{c}{\sqrt{L_0}} \cdot \brr*{K-2 - \abs{\cT_2(K-2)} - \abs{\cT_4(K-2)} - 1}
  \\&=
  \tfrac{c}{\sqrt{L_0}} \cdot \brr*{K+1 - \abs{\cT_2(K-2)} - \abs{\cT_4(K-2)} - 1}
  \\&\at{uses the definition of $\cT_2(k)$ and $\cT_4(k)$ in \cref{eq:Tau}}\geq
  \tfrac{c}{\sqrt{L_0}} \cdot \brr*{K+1 - \abs{\cT_2(K+1)} - \abs{\cT_4(K+1)} - 1},
\end{align*}
where \annotate. This proves Case~2.

\textbf{Case 3a.}
Since $l(K+1) = 0$, we have $\cT_1(K+1) = \cT_2(K+1) = \varnothing$. Hence, from \cref{eq:Tau,eq:l} and \cref{line:eta_H}, we have $\eta_{K+1} = (1+\gamma)^{K+1}\eta_0$. Moreover, since $K+1 \in \cT_3(K+1)$ and $l(K+1) \notin \cT_1(K+1)$, we have $K+1 > m$ due to \cref{eq:Tau,eq:l}. Besides, since $\cT_1(K+1) = \cT_2(K+1) = \varnothing$, from \cref{eq:Tau} we conclude that $\cT_3 (K+1) = \brf{m+1,\dots,K+1}$ and $\cT_4 (K+1) = \brf{1,\dots, m}$. Hence, we can lower-bound $\sqrt{H_{K+1}}$ as follows:
\begin{align*}
  \sqrt{H_{K+1}}
  &\at{is obtained similarly to \cref{aux:L:2a} in the proof of \Cref{thm:L} in \Cref{proof:thm:L}}\geq
  \tfrac{\eta_0\gamma(1+\gamma)^m}{4}\brr*{K-m+2}^2
  \\&\at{uses the definition of $m$ in \cref{eq:cm_L01}}\geq
  \tfrac{c^2}{L_0}\brr*{K-m+2}^2
  \\&\at{uses the fact that $\cT_4 (K+1) = \brf{1,\dots, m}$ and $\cT_2(K+1) = \varnothing$ as shown above}=
  \tfrac{c^2}{L_0}\brr*{K+2-\abs{\cT_2(K+1)}-\abs{\cT_4(K+1)}}^2
  \\&\geq
  \tfrac{c^2}{L_0}\brr*{K+1-\abs{\cT_2(K+1)} - \abs{\cT_4(K+1)} - 1}^2,
\end{align*}
where \annotate. This proves Case~3a.

\textbf{Case 3b.}
Using \cref{eq:Tau,eq:l} and \cref{line:eta_H}, we can express $\eta_{K+1}$ as follows:
\begin{equation}
  \eta_{K+1}
  =
  (1+\gamma)^{K+1 - l(K+1)} \eta_{l(K+1)}
  =
  \tfrac{\nu H_{l(K+1)-2}\lambda_{l(K+1)}}{\eta_{l(K+1)-2}}.
\end{equation}
Moreover, similar to \cref{aux:L:H} in the proof of \Cref{thm:L} in \Cref{proof:thm:L}, we can upper-bound $H_{K+1}$ as follows:
\begin{equation}
  H_{K+1} \leq \tfrac{1}{\gamma}(1+\gamma)^{K-l(K+1) + 3}H_{l(K+1)-1}.
\end{equation}
In addition, similar to \cref{aux:L:2b} in the proof of \Cref{thm:L} in \Cref{proof:thm:L}, and using the definition of $\cT_1(k)$ in \cref{eq:Tau}, we can obtain the following inequality:
\begin{equation}
  \brr*{\sqrt{H_{K+1}} - \sqrt{H_{l(K+1)-3}}}^2 \geq \tfrac{\nu}{18 L_0\sqrt{\gamma(1+\gamma)(2+\gamma)^3}}\cdot\brr*{(1+\gamma)^{(K-l(K+1))/2+1} - 1},
\end{equation}
where $H_{-2} = H_0$, and similar to  \cref{aux:L:2b_1} in the proof of \Cref{thm:L} in \Cref{proof:thm:L}, we can obtain the following inequality:
\begin{equation}\label{aux:L01:3b}
  \sqrt{H_{K+1}} - \sqrt{H_{l(K+1)-3}}
  \geq
  \tfrac{\sqrt{\nu}\gamma}{24\sqrt[4]{4\gamma(1+\gamma)^5(2+\gamma)^3L_0^2}}\cdot\brr*{K-l(K+1)+4}.
\end{equation}
Finally, after rearranging, we can lower-bound $\sqrt{H_{K+1}}$ as follows:
\begin{align*}
  \sqrt{H_{K+1}}
  &\at{uses \cref{aux:L01:3b}}\geq
  \sqrt{H_{l(K+1)-3}} + \tfrac{\sqrt{\nu}\gamma}{24\sqrt[4]{4\gamma(1+\gamma)^5(2+\gamma)^3L_0^2}}\cdot\brr*{K-l(K+1)+4}
  \\&\at{uses the definition of $c$ in \cref{eq:cm_L01}}\geq
  \sqrt{H_{l(K+1)-3}} + \tfrac{c}{\sqrt{L_0}}\cdot\brr*{K-l(K+1)+4}
  \\&\at{uses the induction hypothesis in \cref{eq:H_L01}}\geq
  \tfrac{c}{\sqrt{L_0}} \cdot \brr*{l(K+1)-3 - \abs{\cT_2(l(K+1)-3)} - \abs{\cT_4(l(K+1)-3)} - 1}.
  \\&
  +\tfrac{c}{\sqrt{L_0}}\cdot\brr*{K-l(K+1)+4}
  \\&=
  \tfrac{c}{\sqrt{L_0}}\cdot\brr*{K + 1 - \abs{\cT_2(l(K+1)-3)} - \abs{\cT_4(l(K+1)-3)} - 1}
  \\&\at{uses the definition of $\cT_2(k)$ and $\cT_4(k)$ in \cref{eq:Tau}}\geq
  \tfrac{c}{\sqrt{L_0}}\cdot\brr*{K + 1 - \abs{\cT_2(K+1)} - \abs{\cT_4(K+1)} - 1},
\end{align*}
where \annotate. This proves Case~3b.

\textbf{Case 3c.}
Similar to \cref{aux:L:2b} in the proof of \Cref{thm:L} in \Cref{proof:thm:L}, and using \Cref{lem:lambda_min}, we can obtain the following inequality:
\begin{align*}
  \brr*{\sqrt{H_{K+1}} - \sqrt{H_{l(K+1)-3}}}^2
  &\geq
  \tfrac{\nu\lmin}{4\sqrt{\gamma(1+\gamma)(2+\gamma)^3}}\cdot\brr*{(1+\gamma)^{(K-l(K+1))/2+1} - 1}
  \\&=
  \tfrac{\nu\lmin(1+\gamma)^{m/4}}{4\sqrt{\gamma(1+\gamma)(2+\gamma)^3}}\cdot\brr*{(1+\gamma)^{(K-l(K+1) - m/2)/2+1} - \tfrac{1}{(1+\gamma)^{m/4}}}
  \\&\at{uses the definition of $m$ in \cref{eq:cm_L01}}\geq
  \tfrac{\nu}{18L_0\sqrt{\gamma(1+\gamma)(2+\gamma)^3}}\cdot\brr*{(1+\gamma)^{(K-l(K+1) - m/2)/2+1} - 1}
  \\&\at{uses the fact that $K - l(K+1) \geq m$, which is implied by the assumptions $K+1 \in \cT_3(K+1)$ and $l(K+1)\in\cT_2(K+1)$, and \cref{eq:Tau}}\geq
  \tfrac{\nu}{18L_0\sqrt{\gamma(1+\gamma)(2+\gamma)^3}}\cdot\brr*{(1+\gamma)^{(K-l(K+1))/4+1} - 1},
\end{align*}
where $H_{-2} = H_0$, and \annotate.
After taking the square root from both sides of the inequality, we obtain the following:
\begin{align*}
  \sqrt{H_{K+1}} - \sqrt{H_{l(K+1)-3}}
  &\geq
  \tfrac{\sqrt{\nu}}{3\sqrt[4]{4\gamma(1+\gamma)(2+\gamma)^3L_0^2}}\cdot\sqrt{(1+\gamma)^{(K-l(K+1))/4+1} - 1}
  \\&\at{uses the inequality $\sqrt{a}\geq \sqrt{a+b} - \sqrt{b}$}\geq
  \tfrac{\sqrt{\nu}}{3\sqrt[4]{4\gamma(1+\gamma)(2+\gamma)^3L_0^2}}\cdot\brr*{(1+\gamma)^{(K-l(K+1)+4)/8} - 1}
  \\&=
  \tfrac{\sqrt{\nu}}{3\sqrt[4]{4\gamma(1+\gamma)^5(2+\gamma)^3L_0^2}}\cdot\brr*{(1+\gamma)^{(K-l(K+1)+4)/8+1} - (1+\gamma)}
  \\&\at{uses Bernoulli's inequality}\geq
  \tfrac{\sqrt{\nu}\gamma}{24\sqrt[4]{4\gamma(1+\gamma)^5(2+\gamma)^3L_0^2}}\cdot\brr*{K-l(K+1)+4},
\end{align*}
where \annotate. The rest of the proof is identical to the proof of Case~3b.\qed

\proofsubsection{cor:L01}

We can upper-bound $\cD^2$ as follows:
\begin{align*}
  \cD
  &\at{uses \cref{eq:D}}=
  \sqrt{\max\brf*{\tfrac{1}{\gamma\theta},2(1+\gamma), (1+2\theta)^2}
    \brr*{
      \sqn{x_0 - x^*}
  +(1+\gamma\theta)\eta_0^2\sqn{\nabla f(x_0)}}}
  \\&\at{uses the inequality $\sqrt{a+b}\leq \sqrt{a}+\sqrt{b}$}\leq
  \sqrt{\max\brf*{\tfrac{1}{\gamma\theta},2(1+\gamma), (1+2\theta)^2}}
  \brr*{
    \norm{x_0 - x^*}
    +\sqrt{(1+\gamma\theta)}\eta_0\norm{\nabla f(x_0)}
  }
  \\&\at{uses \cref{eq:L01}}\leq
  \sqrt{\max\brf*{\tfrac{1}{\gamma\theta},2(1+\gamma), (1+2\theta)^2}}
  \brr*{
    1+
    (1+\gamma\theta)\eta_0L_0\cdot\tfrac{\brr{\exp(L_1\norm{x_0-x^*}) - 1}}{L_1\norm{x_0 - x^*}}
  }\norm{x_0 - x^*}
  \\&\at{uses the inequality $ \exp(t)-1 \leq t\exp(t)$, which is implied by the convexity of the function $t \mapsto \exp(t)$}\leq
  \sqrt{\max\brf*{\tfrac{1}{\gamma\theta},2(1+\gamma), (1+2\theta)^2}}
  \brr*{
    1+
    (1+\gamma\theta)\eta_0L_0\exp(L_1\norm{x_0-x^*})
  }\norm{x_0 - x^*}
  \\&\at{uses the assumption $\eta_0L_0 \leq \exp(-L_1\norm{x_0-x^*})$}\leq
  \sqrt{\max\brf*{\tfrac{1}{\gamma\theta},2(1+\gamma), (1+2\theta)^2}}
  \brr*{
    2+\gamma\theta
  }\norm{x_0 - x^*}
  \\&=
  \cO\brr{\norm{x_0-x^*}},
\end{align*}
where \annotate. It remains to combine \Cref{cor:main}, \Cref{lem:Tau}, and \Cref{thm:L01}.\qed

% \proofsubsection{rem:L01}

% We can lower-bound $\norm{x_0 - x^*} / \norm{\nabla f(x_0)}$ as follows:
% \begin{align*}
%   \norm{x_0 - x^*} / \norm{\nabla f(x_0)}
%   &\at{uses \cref{eq:L01}}\geq
%   \tfrac{1}{L_0} \cdot \tfrac{L_1\norm{x_0-x^*}}{\exp(L_1\norm{x_0-x^*}) - 1}
%   \\&\at{uses the inequality $t\exp(t) \geq \exp(t)-1$, which is implied by the convexity of the function $t \mapsto \exp(t)$}\geq
%   \tfrac{1}{L_0}\exp(-L_1\norm{x_0-x^*})
%   \\&= \eta_0,
% \end{align*}
% where \annotate. It remains to plug the stepsize $\eta_0 = \tfrac{1}{L_0}\exp(-L_1\norm{x_0-x^*})$ into \cref{eq:K_L01}.\qed

% \input{l01.tex}

\end{document}